\documentclass[12pt]{article}
\usepackage[fleqn]{amsmath}
\usepackage[cp1251]{inputenc}
\usepackage[russian,english]{babel}
\usepackage{amsmath,amsfonts,amssymb,graphicx,makeidx,amsthm,cite}
\usepackage{xcolor}
\usepackage{hyperref}
\usepackage[mathscr]{euscript}
\let\euscr\mathscr \let\mathscr\relax
\usepackage[scr]{rsfso}
\hypersetup{
	colorlinks = false,
	linkbordercolor = {violet},
}

\newtheorem{rem}{Remark}
\newtheorem{criteria}{Criteria}

\begin{document}
 
\title{{\bf Overview of fractional calculus and its computer implementation in Wolfram Mathematica}}

\author{{\bf  Oleg Marichev} \\(Wolfram Research, 100 Trade Center Drive
        Champaign,\\ IL 61820-7237, USA, oleg@wolfram.com) 
\and
        {\bf  Elina Shishkina} \\(Voronezh State University,\\ Universitetskaya Pl., 1, 394018, Voronezh,  Russia,\\
       Belgorod State National Research University,\\ Pobeda Street, 85, 308015, Belgorod, Russia, shishkina@amm.vsu.ru)
 }

\maketitle

\thispagestyle{empty}
\newpage

\begin{abstract}
	This survey  aims to present  various approaches to non-integer integrals and derivatives and  their practical implementation within Wolfram Mathematica. It begins by  short discussion of historical moments and applications related to fractional calculus. Different methods for handling non-integer powers of differentiation operators are presented, along with generalizations of fractional integrals and derivatives. The survey also delves into the diverse applications of fractional calculus in physics, engineering, medicine, and numerical calculations.
	
	Essential details of fractional integro-differentiation implemented in Wolfram Mathematica are highlighted. The Hadamard regularization of  Riemann-Liouville operator is utilized as the foundation for creating the arbitrary order of integro-differential operator in Mathematica. The survey describes the application of fractional integro-differentiation to Taylor series expansions near zero using Hadamard regularization and the use of the Meijer $G$-function for evaluating derivatives of complex orders.
	
	We conclude with a discussion on applying fractional integro-differentiation to "differential constants"\, and provide generic formulas for fractional differentiation. The extensive list of references underscores the vast body of works on fractional calculus.

	
	
\end{abstract} 




\setcounter{section}{0} \setcounter{equation}{0} 

\section{Introduction} \label{sec:1}

In this survey, we would like to consider some approaches to non-integer integro-differentiation and its implementation in the computer algebra system Wolfram Mathematica.

In Section \ref{S01}, some historical moments and applications are briefly discussed. The ideas that initiated the era of fractional calculus are presented in Subsection \ref{S01.1}. Different approaches to non-integer powers of $\frac{d}{dx}$ are provided in Subsection \ref{DA}. Generalizations of fractional derivatives, such as operators with a power-logarithmic kernel, non-integer powers of  $x\frac{d}{dx}$, and $\frac{1}{x^\gamma}\frac{d}{dx}x^\gamma\frac{d}{dx}$, fractional derivatives of a function $f$ with respect to another function $g$, Erd\'{e}lyi-Kober type operators, variable-order differential operators, and sequential fractional derivatives are described in Subsection \ref{GFD}. Fractional calculus has diverse and widespread applications. References to some of numerous works  where fractional calculus is applied to physics, engineering, medicine, as well as some approaches to the numerical calculation of fractional derivatives are presented in Subsection \ref{Apl}.

In Section \ref{WS}, some essential details of fractional integro-differentiation, which have been implemented in Wolfram Mathematica, are presented. This section forms a new perspective on fractional differentiation for analytic functions of complex variables. The Hadamard regularized Riemann-Liouville operator serves as the basis for creating the integro-differential operator in Wolfram Mathematica. The Riemann-Liouville-Hadamard integro-differentiation of an arbitrary function to an arbitrary symbolic order $\alpha$ in the Wolfram Language is described in Subsection \ref{Had}. Subsection \ref{HR} presents an approach in which fractional integro-differentiation is applied to each term of Taylor series expansions of all functions near zero using Hadamard regularization. For the evaluation of the derivative of an arbitrary complex order, the Meijer $G$-function can be utilized. Many functions can be represented as a Meijer $G$-function, facilitating the finding of an $\alpha$-derivative of the Meijer $G$-function. This approach is detailed in Subsection \ref{MG}. Subsection \ref{BF} discusses how fractional integro-differentiation is applied to "differential constants", which are not constants everywhere but vary on some domains of the complex plane, exhibiting discontinuity on certain lines. This subsection also includes the description of generic formulas for fractional differentiation.

The list of references contains numerous items, but the total number of works on fractional calculus   is however immense, indeed.


\section{ History, nowadays and applications}\label{S01}

\subsection{Background and concept}\label{S01.1}

\subsubsection{Gamma function and the idea of non-integer derivative}

Calculus teaches us how to compute derivatives of any integer order. We can interpret the differentiation of negative integer order as repeated integration. The zeroth order of differentiation gives the function itself. The question is how to generalize derivatives to non-integer orders.

For extensive details on the history of fractional calculus, we refer readers to the handbook \cite{SamkKilbMari1993}, especially to the numerous historical notes there. More information on the main milestones of fractional calculus can be found in \cite{LazarevHist,Ross}. The history of fractional calculus in recent times is presented in \cite{MKM,MK1} and in \cite{Handbook}, Vol. 1, 1--22 (2019). The first application of fractional calculus was discovered by Abel (see \cite{PodlubnyAbel} for details).

The starting point for the transition from integer-order to fractional-order derivatives and integrals is Euler's introduction of the Gamma function. Namely, in 1729, Euler derived the integral representation for the factorial $n!$, leading to the Gamma function $\Gamma(n+1)=n!$, allowing the definition of $n!$ for any complex number $n$ except $0, -1, -2, \ldots$. Besides the Gamma function, there are other ways to generalize the factorial to complex numbers (see \cite{Marichev}, p. 35-36, \cite{Marichev01}).

Euler introduced the Gamma function through the integral, which later became known as the Euler integral of the second kind:
\begin{equation}\label{Gamma}
	\Gamma (\alpha)=\int\limits_{0}^{\infty }x^{\alpha-1}e^{-x}\,dx,\qquad {\rm Re}(\alpha)>0.
\end{equation}
The property $\Gamma (\alpha)=\frac {\Gamma (\alpha+1)}{\alpha}$ can be used to uniquely extend $\Gamma (\alpha)$ to a meromorphic function defined for all complex numbers $\alpha$ except integers less than or equal to zero. If $\alpha{=}0,{-}1,{-}2,{...}$ the Gamma function has simple poles. 

The Bohr--Mollerup theorem \cite{Bohr} shows that the Gamma function is the only function that satisfies the properties
\begin{enumerate}
	\item $f(1)=1$,
	\item $f(x+1)=xf(x)$,
	\item for every $x\geq 0$
	$\ln f$	is a convex function.
\end{enumerate}

So  the Euler's Gamma function is the "best"\, extension of the factorial function to the real (complex) numbers.
The Gamma function in Wolfram Mathematica, denoted as \texttt{Gamma[z]}, is suitable for both symbolic and numerical manipulation, $z\in\mathbb{C}$.

This generalization of the factorial allowed Euler to realize that the concept of the $n$-th order derivative of the power function $x^p$ acquired meaning for a non-integer $n$.

Namely, let $n\in\mathbb{N}$, $x>0$ and $p\in\mathbb{R}$. It is obvious that
$$
(x^p)^{(n)}=p(p-1)...(p-n+1)x^{p-n}
$$
or
\begin{equation}\label{eq1}
	(x^p)^{(n)}=\frac{p!}{(p-n)!}x^{p-n}.
\end{equation}
The expression \eqref{eq1} has sense for non-integer $n$ when we use the Gamma function \eqref{Gamma}. The derivative of $x^p$ of non-integer order $\alpha$ can be defined by
\begin{equation}\label{eq2}
	\frac{d^\alpha}{dx^\alpha}\,x^p=\frac{\Gamma(p+1)}{\Gamma(p-\alpha+1)}\,x^{p-\alpha} .
\end{equation}

It should be noted that formula \eqref{eq2} is valid  not for all values $p$ and $\alpha$. Therefore, we need to exclude the arguments that result in poles of the Gamma function in the numerator. That means that in \eqref{eq2} 
$$-(p+1)\notin\mathbb{N}\cup\{0\}.$$

We know that for $\alpha=0$ we just get the $\frac{d^0}{dx^0}\,x^p=x^p$ and for $\alpha=-1$ we get integral and for $\alpha=-2,-3,...$ iterated integrals.  For example, if $p=-1$ and $\alpha=-1$ we obtain 
$$
\frac{d^{-1}}{dx^{-1}}\,x^{-1}=\int\frac{dx}{x}=\log(x)+C.
$$
This case already results in a logarithm rather than a power function.  Let us discuss how we get a logarithm here.
We can say that this is obtained by regularizing the integral $\int x^pdx$ and then taking a limit $p\rightarrow -1$.

Let us see how to obtain the logarithm from formula below, including an arbitrary constant $C$:
$$\int x^pdx=\frac{x^{p+1}}{p+1}+C$$
by taking the limit. 

We can change this constant to another constant $-\frac{1}{p+1}{+}C_1$, including parameter 
$p$ and a constant $C_1$.

After evaluating the integral and taking the limit, we get
$$
\int x^pdx=\frac{x^{p+1}}{p+1}-\frac{1}{p+1}+C_1=\frac{x^{p+1}-1}{p+1}+C_1,
$$
and also
$$
\int\frac{dx}{x}=\lim\limits_{p\rightarrow-1}\int x^pdx=\lim\limits_{p\rightarrow-1}\frac{x^{p+1}-1}{p+1}+C_1=\log(x)+C_1.
$$
Here, we applied regularization using the term $\left( -\frac{1}{p+1}\right) $. 
Integration implemented in Wolfram Mathematica also uses regularization when necessary.
The general approach to finding an arbitrary order derivative of a power function is considered in Subsection \ref{HR}.


Therefore, if function $f(x)$ is locally given by a convergent power series or $f(x)$ is an analytic function:
$$
f(x)=\sum\limits_{p=0}^{\infty }a_px^p,\qquad a_p={\frac {f^{(p)}(0)}{p!}},
$$
then the derivative of order $\alpha>0$ can be formally defined as
$$
\frac{d^\alpha f(x)}{dx^\alpha}=\sum\limits_{p=0}^{\infty }\frac{\Gamma(p+1)}{\Gamma(p-\alpha+1)}\,a_px^{p-\alpha}.
$$
To ensure the existence of such a fractional derivative, it is necessary to verify the convergence of the specified series. Here,   $p\geq 0$ allow to avoid the poles of the Gamma function in the numerator. The general case when $\alpha$ and $p$ can be negative is considered in Subsection \ref{HR}.

In 1823, the first application of fractional calculus was discovered by Abel (see for details \cite{PodlubnyAbel,LazarevHist}). Abel was searching for a curve in the plane such that the time required for a particle to slide down the curve to its lowest point under the influence of gravity is independent of its initial position on the curve. Such a curve is called the tautochrone. Abel found that to find this curve, it is necessary to solve the equation
$$
\int\limits_0^x\frac{f(t)dt}{(x-t)^{1/2}}=\varphi(x)= T ={\rm const},\qquad x>0.
$$
It is worth noting that Abel solved a more general equation
$$
\int\limits_0^x\frac{f(t)dt}{(x-t)^{\alpha}}=\varphi(x), \qquad x>0,
$$
where $0<\alpha<1$.

Furthermore, in 1832, Liouville formally extended formula for the integer derivative of the exponent $\frac{d^n}{dx^n}e^{bx}$ ($b$ is some number) to derivatives of arbitrary order $\frac{d^\alpha}{dx^\alpha}e^{bx}$. Namely,
\begin{equation}\label{Leibniz}
	\frac{d^\alpha e^{bx}}{dx^\alpha}=b^\alpha e^{bx}.
\end{equation}
Based on formula \eqref{Leibniz}, one can formally write the derivative of order $\alpha\in\mathbb{R}$ of an arbitrary function $f$ represented by the series
$$
\frac{d^\alpha f(x)}{dx^\alpha}=\sum\limits_{k=0}^\infty c_k b_k^\alpha e^{b_k x},\qquad \text{where}\qquad f(x)=\sum\limits_{k=0}^\infty c_k e^{b_k x}.
$$
The limitation of this definition is related to the convergence of the series.

\subsubsection{Formal definition of fractional integro-differentiation}

So, starting from the 17th century, the need for a formal definition of fractional integro-differentiation gradually formed. We present here such a definition following \cite{Ross}.

Let $z\in\mathbb{C}$, $x\in\mathbb{R}$  be variables. 
The starting point in fractional integro-differentiation is $0$, $\nu\in\mathbb{C}$ or $\nu\in\mathbb{R}$ is an order of integro-differentiation.
Fractional integro-differentiation $D^\nu$ acts to $f(z)$ by $z$ and acts to $f(x)$ by $x$. 
In \cite{Ross}, one can find the following criteria of fractional integro-differentiation.

\begin{criteria}\label{cr01} Operator $D^\nu f(z)$ is the integro-differential operator of order $\nu\in\mathbb{C}$ if and only if
	\begin{enumerate}
		\item If $f(z)$ is an analytic function of the complex variable
		$z$, the derivative $D^\nu f(z)$ is an analytic function of $\nu$ and $z$.
		\item The operation $D^\nu f(z)$ must produce the same result
		as ordinary differentiation when $\nu$ is a positive integer. If $\nu$
		is a negative integer, say $\nu=-n$, then $D^{-n} f(z)$
		must produce the 	same result as ordinary $n$-fold integration and $D^{-n} f(z)$ must vanish 
		along with its $(n-1)$-derivatives at $z=0$.
		\item The operation of order zero leaves the function unchanged:
		$$
		D^{0} f(z)=f(z).
		$$
		\item The fractional operators must be linear:
		$$
		D^\nu [af(z)+bg(z)]= aD^\nu f(z)+bD^\nu g(z).
		$$
		\item The  semigroup  property of arbitrary order for suitable function $f$ holds:
		$$
		D^\nu\,	D^\mu f(z)=	D^{\nu+\mu} f(z).
		$$
	\end{enumerate}
\end{criteria}

For fractional integral operators, the semigroup property is most often satisfied for all functions for which this integral exists. However, for fractional derivatives, it is necessary to construct a suitable class of functions to satisfy the semigroup property (see, for example, \cite{SamkKilbMari1993}, p. 34 and p. 45 and further).

In the next Subsections \ref{DA} and \ref{GFD}, we examine different operators that realize non-integer powers of differentiations and integrations.

\subsection{Different approaches to classical fractional calculus}\label{DA}

In this subsection, we consider different classical approaches to introduce a non-integer power $\left(\frac{d}{dx} \right)^\alpha$
of the differentiation operator$\frac{d}{dx}$. If an operator has two forms, left-sided and right-sided, we will consider only the left-sided.

\subsubsection{Riemann-Liouville type fractional integrals and derivatives}

The most well known definitions of operators  giving  arbitrary real (complex) powers of the $\frac{d}{dx}$
which fulfill  Criteria \ref{cr01} are
\begin{equation}\label{Eq01}
	(I_{a+}^{\alpha } f )(x)=\frac{1}{\Gamma (\alpha )} \int\limits_{a}^{x}\frac{f(t)}{(x-t)^{1-\alpha}}dt ,\quad a<x\leq b,\quad \alpha>0,	
\end{equation}
and
\begin{equation} \label{RLD3} 
	(D_{a+}^{\alpha }f)(x)=\frac{1}{\Gamma (n{-}\alpha )} \left(\frac{d}{dx} \right)^{n} \int\limits_{a}^{x}\frac{f(t)dt}{(x-t)^{\alpha -n+1} } ,\quad a<x\leq b,\quad \alpha>0,	
\end{equation}
where $n=[\alpha]+1$, when $\alpha$ in not integer and $n=\alpha$, when $\alpha\in\mathbb{Z}_+$.
The definition of the operator \eqref{Eq01} is valid for ${\rm Re}(\alpha) >0$ and the definition of the operator \eqref{RLD3} is valid for $n-1 \le {\rm Re}(\alpha) < n$.

Operators \eqref{Eq01} and \eqref{RLD3} are called left-sided Riemann-Liouville integral and derivative on the segment $[a,b]$, respectively. 
The definition  of fractional integral is based on a generalization of the formula for an $n$-fold integral 
\begin{equation}\label{Cauchy}
	\underbrace{\int\limits_a^xdx...\int\limits_a^{x}dx\int\limits_a^{x}}_nf(x)dx=\frac{1}{(n-1)!}\int\limits_a^x(x-t)^{n-1}f(t)dt ,\qquad x>a.
\end{equation}

Then the fractional derivative \eqref{RLD3} is obtained as a left inverse operator to \eqref{Eq01} which is constructed as a solution to the Abel equation. It is necessary to note that the operator \eqref{RLD3} is one of infinitely many different left-inverse operators to the fractional Riemann-Liouville integral  \eqref{Eq01}.

If we consider $D_{a+}^{\alpha } x^p$, $x>0$, $p\in\mathbb{R}$ we obtain
\begin{align*}
	D_{a+}^{\alpha } x^p&=\frac{1}{\Gamma (n-\alpha )} \left(\frac{d}{dx} \right)^{n} \left( \int\limits_{a}^{x}\frac{t^pdt}{(x-t)^{\alpha -n+1} }\right)\\
	&=\frac{1}{\Gamma (n-\alpha )} \left(\frac{d}{dx} \right)^{n} x^{p-\alpha+n} \left(\frac{\Gamma (p+1) \Gamma (n-\alpha )}{\Gamma (p+n-\alpha +1)}{-}B_{\frac{a}{x}}(p+1,n-\alpha )\right)\\
	&=\frac{\Gamma(p+1)}{\Gamma(p-\alpha+1)}x^{p-\alpha}-\frac{1}{\Gamma (n-\alpha )} \left(\frac{d}{dx}\right)^{n} x^{p-\alpha+n} B_{\frac{a}{x}}(p+1,n-\alpha ),
\end{align*}
where $B_z(a,b)$ is  the incomplete Beta function 
$$
B_z(a,b)=\int\limits_0^z t^{a-1}(1-t)^{b-1}dt.
$$

It is easy to see that in order for the resulting formula to correspond to \eqref{eq2}, we need to set $a=0$,  when the last term  disappears: 
$$D_{0+}^{\alpha } x^p=\frac{\Gamma(p+1)}{\Gamma(p-\alpha+1)}x^{p-\alpha}.$$
Because of this, in   Section \ref{WS} we will consider only fractional derivatives and integrals tied to the origin, i.e. with $a=0$.

Liouville fractional integral of order $\alpha$ on the semiaxis $(0,\infty)$ for $f{\in}L_1(0,\infty)$ has the form
\begin{equation}\label{Eq03}
	(I_{0+}^{\alpha } f )(x)=\frac{1}{\Gamma (\alpha )} \int\limits_{0}^{x}\frac{f(t)}{(x-t)^{1-\alpha}}dt,\qquad x\in(0,\infty).
\end{equation}

Liouville fractional integral of order $\alpha$ on  the whole real axis is 
\begin{equation}\label{Eq05}
	(I_{+}^{\alpha } f )(x)=\frac{1}{\Gamma (\alpha )} \int\limits_{-\infty}^{x}\frac{f(t)}{(x-t)^{1-\alpha}}dt ,\qquad x\in \mathbb{R}.
\end{equation}

Liouville fractional derivatives of order $\alpha$ on  the semi-axis and on the  whole real axis  are, respectively 
\begin{equation} \label{RLD5} 
	(\mathcal{D}_{0+}^{\alpha }f)(x)=\frac{1}{\Gamma (n-\alpha )} \left(\frac{d}{dx} \right)^{n} \int\limits_{0}^{x}\frac{f(t)dt}{(x-t)^{\alpha -n+1} },   \qquad x\in(0,+\infty),
\end{equation}
\begin{equation} \label{RLD7} 
	(\mathcal{D}_{+}^{\alpha }f)(x)=\frac{1}{\Gamma (n-\alpha )} \left(\frac{d}{dx} \right)^{n} \int\limits_{-\infty}^{x}\frac{f(t)dt}{(x-t)^{\alpha -n+1} }, \qquad x\in \mathbb{R},   
\end{equation}
where $n=[\alpha ]+1$, when $\alpha$ in not integer and $n=\alpha$, when $\alpha\in\mathbb{Z}_+$.

Next, we consider a  modification of the Riemann-Liouville integration, for which lower and upper integral limits are symmetric. 
For fixed arbitrary point $c\in\mathbb{R}$ and $\alpha>0$ operator
\begin{equation} \label{Chen} 
	(I_c^\alpha f)(x)=
	\frac{1}{\Gamma (\alpha )}\left\{ \begin{array}{ll}
		\int\limits_{c}^{x}(x-t)^{\alpha-1}f(t)dt & \mbox{if $x>c$};\\
		\int\limits_{x}^{c}(t-x)^{\alpha-1}f(t)dt & \mbox{if $x<c$}\end{array} \right.
\end{equation}
is called the Chen fractional integral.

Let $\alpha>0$ be not integer, $n=[\alpha ]+1$ then Chen fractional derivative is 
\begin{equation} \label{ChenDer01} 
	(D_c^\alpha f)(x)=
	\frac{1}{\Gamma (n-\alpha )}\left\{ \begin{array}{ll}
		\left(\frac{d}{dx} \right)^{n} \int\limits_{c}^{x}\frac{f(t)dt}{(x-t)^{\alpha -n+1} } & \mbox{if $x>c$};\\
		\left(-\frac{d}{dx} \right)^{n} \int\limits_{x}^{c}\frac{f(t)dt}{(t-x)^{\alpha -n+1} } & \mbox{if $x<c$.}\end{array} \right.
\end{equation}

One of the most popular fractional derivative that is used now is the Caputo fractional derivative 
which has the form
\begin{equation}\label{Cap0}
	(\,^{C}D^\alpha_{0+} f)(x)=\frac{1}{\Gamma(n-\alpha)}\int\limits_{0}^x \frac{f^{(n)}(y)dy}{(x-y)^{\alpha-n+1}},\qquad n=[\alpha]+1,\qquad x>0.
\end{equation}
Operator \eqref{Cap0}
with $n=1$ was introduced already by Abel, see the paper \cite{PodlubnyAbel}.
Original works of   Italian  physicist  M. Caputo  are \cite{Caputo1} and   \cite{Caputo2}.   

In 1948 (see \cite{Gerasimov}, submitted in 1947) the Soviet  mechanics A.N. Gerasimov introduced fractional derivative of the form
\begin{equation}\label{Ger0}
	\frac{1}{\Gamma(1-\alpha)}\int\limits_{-\infty}^x
	\frac{f'(y)dy}{(x-y)^\alpha},\quad y>0,\quad x\in\mathbb{R},
	\quad 0<\alpha<1.
\end{equation}
In the same work, A.N. Gerasimov studied two new problems in
viscoelasticity theory. He reduced this problem to  differential equations with partial fractional derivative.

Let us notice that in the article \cite{PodlubnyAbel} Abel's contribution to the theory of fractional calculus was carefully analyzed, in particular the appearance of the derivative \eqref{Cap0}  for $\alpha\in(0,1)$.

The popularity of the fractional derivative  \eqref{Cap0} in applications is explained as follows.
If we consider the fractional differential equation with Riemann--Liouville fractional derivative
of the form
$$
(D^{\alpha}_{0+}f)(x)=\lambda f(x),\qquad x>0,\qquad 0<\alpha< 1,\qquad \lambda\in\mathbb{R}
$$
we should add the initial condition
$$
(D^{\alpha-1}_{0+}f)(0+)=1
$$
and the solution of this problem is (see \cite{KilSrTr})
\begin{equation}\label{RLCP1}
	f(x)= x^{\alpha-1}E_{\alpha,\alpha}(\lambda x^\alpha)
\end{equation}
where
$
E_{\alpha,\beta}(\lambda x^\alpha)=\sum\limits_{k=0}^{\infty}\frac{(\lambda x^\alpha)^k}{\Gamma(\alpha k+\beta)}
$
is the  Mittag--Leffler function (see \cite{Gorenflo}, \cite{Handbook}, Vol. 1, 269--296 (2019)). Note that we have a singularity at zero in \eqref{RLCP1}. 
Therefore, it is not possible to consider the classical 
Cauchy problem at zero initial point with Riemann--Liouville fractional derivative.

From the other side, the   solution to the Cauchy problem for fractional differential equation with  fractional derivative \eqref{Cap0}
$$
(\,^{C}D^{\alpha}_{0+}f)(x)=\lambda f(x),\qquad x>0,\qquad 0<\alpha\leq 1,\qquad \lambda\in\mathbb{R}
$$
$$
f(0+)=1,
$$
is  (see \cite{KilSrTr})
\begin{equation}\label{ML01}
	f(x)= E_{\alpha,1}(\lambda x^\alpha).
\end{equation}
Now it is bounded at zero and so the classical Cauchy problem
is correct.

Let us find  the Riemann--Liouville derivative of \eqref{ML01}:
\begin{align*}
	D^{\alpha}_{0+}E_{\alpha,1}(\lambda x^\alpha)&=\sum_{k=0}^{\infty}\frac{\lambda^k}{\Gamma(\alpha k+1)}D^{\alpha}_{0+} 
	x^{\alpha k}\\
	&= \sum_{k=0}^{\infty}\frac{\lambda^k}{\Gamma(\alpha k+1)}  \frac{\Gamma(\alpha k+1)}{\Gamma(\alpha k-\alpha+1)}x^{\alpha k-\alpha}\\
	&=\sum_{k=0}^{\infty}\frac{\lambda^k x^{\alpha (k-1)}}{\Gamma(\alpha k-\alpha+1)}\\
	&=\frac{x^{-\alpha}}{\Gamma(1-\alpha)}+\lambda \sum_{k=1}^{\infty}\frac{(\lambda x^{\alpha})^{k-1}}{\Gamma(\alpha (k-1)+1)}\\
	& =\frac{x^{-\alpha}}{\Gamma(1-\alpha)}+\lambda\sum_{k=0}^{\infty}\frac{(\lambda x^{\alpha})^{k}}{\Gamma(\alpha k+1)}\\
	&=  \frac{x^{-\alpha}}{\Gamma(1-\alpha)}+\lambda E_{\alpha,1}(\lambda x^\alpha).
\end{align*}
So we see that the term $\frac{x^{-\alpha}}{\Gamma(1-\alpha)}$ has a singularity at $x=0$. 

Then, consider $\,^{C}D^{\alpha}_{0+}x^p$, $p\in\mathbb{R}$:
\begin{align*}
	\,^{C}D^\alpha_{0+} x^p&=\frac{\Gamma(p+1)}{\Gamma(p-n+1)\Gamma(n-\alpha)}
	\int\limits_{0}^x \frac{y^{p-n}dy}{(x-y)^{\alpha-n+1}}\\
	&=\frac{\Gamma(p+1)}{\Gamma(p-n+1)\Gamma(n-\alpha)}\frac{\Gamma (n-\alpha ) \Gamma (p-n+1) x^{p-\alpha }}{\Gamma (p-\alpha +1)}\\
	&=\frac{\Gamma(p+1)}{\Gamma(p-\alpha+1)}x^{p-\alpha}.
\end{align*}
Let us note that integration can only be done with $p>n-1$, where $n=[\alpha]+1$.  Therefore, we have serious limitations and cannot, generally speaking, apply the derivative \eqref{Cap0} to term by term to each member of the power series.
In order to get around this problem, we put
$$
\,^{C}D^\alpha_{0+} x^p= \left\{ \begin{array}{ll}
	\frac{\Gamma(p+1)}{\Gamma(p-\alpha+1)}x^{p-\alpha} & \mbox{if $p>n-1$};\\
	0 & \mbox{if $p\leq n-1$}\end{array} \right. 
$$
then we can find $\,^{C}D^{\alpha}_{0+}E_{\alpha,1}(\lambda x^\alpha)$ by the term-by-term differentiation
of the series
\begin{align*}
	\,^{C}D^\alpha_{0+}E_{\alpha,1}(\lambda x^\alpha)&=\sum_{k=0}^{\infty}\frac{\lambda^k}{\Gamma(\alpha k+1)} \,^{C}D^{\alpha}_{0+} 
	x^{\alpha k}\\
	&=\lambda \sum_{k=1}^{\infty}\frac{(\lambda x^{\alpha})^{k-1}}{\Gamma(\alpha (k-1)+1)}=\lambda E_{\alpha,1}(\lambda x^\alpha).
\end{align*}
So here we saw significant limitations when using a fractional derivative \eqref{Cap0}.

In addition to the Caputo derivative, various other modifications of the Riemann-Liouville  derivative appeared. 
We present here some of them.

Although we only consider left-sided operators, we will mention one right-sided fractional integral of the form 
$$
(I^\alpha_{-}f)(x)=\frac{1}{\Gamma(\alpha)}\int\limits_x^\infty (y-x)^{\alpha-1}f(y)dy.
$$
Some authors call it  Weyl fractional order integral of function $f(y)$ (see \cite{Zhmakin}).


Cossar in \cite{Cossar} presented the following fractional derivative
$$
(\mathcal{D}^\alpha_{-}f)(x)=-\frac{1}{\Gamma(1-\alpha)}
\lim\limits_{N\rightarrow\infty}\frac{d}{dx}\int\limits_x^{N}\frac{f(y)}{(y-x)^{\alpha}}\,dy,\qquad 0<\alpha<1.
$$

Osler fractional derivative (see \cite{Osler}, formula 2.2)
$$
(\,^OD^\alpha_{0+}f)(x)=\frac{\Gamma(1+\alpha)}{2\pi i}\int\limits_{C[a,z^+]}\frac{f(t)}{(t-z)^{\alpha+1}}\,dt,\qquad 0<\alpha<1.
$$
where the contour $C[a,z^+]$  starts and ends at $t=0$.

\subsubsection{Finite differences approach to fractional calculus}

Here we consider the finite differences approach to fractional calculus. This is an approach based on classical limit definition of derivative.

Let us consider how to generalize the derivative of order $n$ in the form
\begin{equation} \label{DC} 
	f^{(n)}(x)=\lim\limits_{h\rightarrow 0}\frac{(\Delta_h^n f)(x)}{h^n},
\end{equation}
where $(\Delta_h^n f)(x)$ is a finite difference 
\begin{equation} \label{FD1} 
	(\Delta_{h}^n f)(x)=\sum\limits_{k=0}^n (-1)^k {\binom {n}{k}}f(x-kh),\qquad {\binom {n}{k}}={\frac {n!}{k!(n-k)!}}
\end{equation}
to the non-integer order.

Since  the Gamma function is a generalization of the factorial function to  non-integer values then for $\alpha>0$ a fractional derivative can be given as a generalization of \eqref{DC} by 
\begin{equation} \label{GLDD} 
	f^{(\alpha)}(x)=\lim\limits_{h\rightarrow +0}\frac{(\Delta_h^\alpha f)(x)}{h^\alpha},
\end{equation}
where
\begin{equation} \label{ID1} 
	(\Delta_{h}^\alpha f)(x)=\sum\limits_{k=0}^\infty (-1)^k {\binom {\alpha}{k}}f(x-kh),\quad \binom {\alpha}{k}=
	\frac{\Gamma(\alpha +1)}{\Gamma(k+1)\Gamma(\alpha - k +1)}.
\end{equation}
For $h>0$ and $\alpha>0$ \eqref{GLDD} is the left-hand sided Gr\"unwald--Letnikov derivative.
We can also consider the case when $h<0$ in \eqref{ID1}.
The left-hand sided Gr\"unwald--Letnikov  fractional 
integral in a real line is
\begin{equation} \label{GLII} 
	f^{(-\alpha)}(x)=\lim\limits_{h\rightarrow +0}h^\alpha(\Delta_h^{-\alpha}f)(x).
\end{equation}

Let us find Gr\"unwald--Letnikov derivative of the power function $x^p$, $x>0$, $p\in\mathbb{R}$.
Since the definition \eqref{GLDD} was given for whole real line, we consider a function
$$
f_p(x)=\begin{cases}x^p,&x>0;\\
	0,&x\leqslant 0.\end{cases}
$$
Then $f_p$ can be written as an inverse Laplace transform
$$
f_p(x)=\frac{\Gamma(p+1)}{2\pi i}\int\limits_{c-i\infty}^{c+i\infty}s^{-p-1}e^{xs}\,ds\qquad(c>0),
$$
then
\begin{align*}
	f_p^{(\alpha)}(x)&=\frac{\Gamma(p+1)}{2\pi i}\lim\limits_{h\rightarrow +0}\frac{1}{h^\alpha}(\Delta_h^\alpha)_x\int\limits_{c-i\infty}^{c+i\infty}s^{-p-1}e^{xs}\,ds\\
	&=\frac{\Gamma(p+1)}{2\pi i}\int\limits_{c-i\infty}^{c+i\infty}e^{xs}s^{-p-1}\lim\limits_{h\rightarrow +0} h^{-\alpha}\sum\limits_{k=0}^\infty (-1)^k {\binom {\alpha}{k}}e^{-khs}ds\\
	&=\frac{\Gamma(p+1)}{2\pi i}\int\limits_{c-i\infty}^{c+i\infty}e^{xs}s^{-p-1}\lim\limits_{h\rightarrow +0} h^{-\alpha}(1-e^{-hs})^\alpha ds\\
	&=\frac{\Gamma(p+1)}{2\pi i}\int\limits_{c-i\infty}^{c+i\infty}e^{xs}s^{\alpha-p-1}ds\\
	&=\frac{\Gamma(p+1)}{\Gamma (p-\alpha +1)}x^{p-\alpha }
\end{align*}
that coincides with $D_{0+}^{\alpha } x^p$ for $x>0$.
Here we consider first $\int\limits_{c-iA}^{c+iA}$, then change limits, then take $A\rightarrow \infty$. Also we used  Lebesgue's dominated convergence theorem.

The definition of Marchaud derivative based on finite differences approach \eqref{ID1} is
\begin{equation} \label{Mar05} 
	(\mathbf{D}_{+}^{\alpha }f)(x)=-\frac{1}{\Gamma (-\alpha)A_l(\alpha)}  
	\int\limits_{0}^{\infty}\frac{(\Delta_{t}^l f)(x)}{t^{\alpha+1} }dt,\quad A_l(\alpha)=\sum\limits_{k=0}^l (-1)^{k-1} {\binom {l}{k}}k^\alpha,
\end{equation}
where $0<{\rm Re}\,\alpha<l$, $l\in\mathbb{N}$ or $l=[\alpha]+1$ when $\alpha\in\mathbb{R}$.

Other modifications to fractional derivatives constructed from fractional differences are presented in Chapter 4 of   book \cite{SamkKilbMari1993}.
For example, Gr\"{u}nwald--Letnikov--Riesz fractional derivative of order $\alpha>0$ was defined by
\begin{equation}\label{GLRD01}
	(\,^{GLR}D_x^\alpha f)(x)=\frac{1}{2\cos(\alpha\pi/2)}\lim\limits_{h\rightarrow 0+}\frac{	(\Delta_{h}^\alpha f)(x)+	(\Delta_{-h}^\alpha f)(x)}{|h|^\alpha}
\end{equation}
where the difference $	(\Delta_{h}^\alpha f)(x)$ of a fractional order $\alpha>0$ is defined by the   series \eqref{ID1}.

Gr\"unwald--Letnikov  operators \eqref{GLDD}  and \eqref{GLII},  Marchaud derivative \eqref{Mar05} and  Gr\"{u}nwald--Letnikov--Riesz fractional derivative \eqref{GLRD01} are considered in \cite{Letnikov,SamkKilbMari1993,Ortigueira01}.

\subsubsection{One-dimensional Riesz and Bessel potentials and their inversions}

Now we consider fractional operators defined as the convolution of a function with the Riesz and Bessel kernels and their inversions.
Such concept comes from the field of potential theory and harmonic analysis.

The integral 
\begin{equation} \label{Riesz01} 
	(I^\alpha f)(x)=\frac{1}{2\Gamma(\alpha)\cos\left(\frac{\pi}{2}\alpha \right) }\int\limits_{-\infty}^\infty \frac{f(t)dt}{|t-x|^{1-\alpha}},\qquad {\rm Re}\,\alpha>0,\qquad \alpha\neq 1,3,5,...
\end{equation}
is called the Riesz potential. 
Along with \eqref{Riesz01}, we consider its modification of the form
\begin{equation} \label{Riesz02} 
	H^\alpha f(x)=\frac{1}{2\Gamma(\alpha)\sin\left(\frac{\pi}{2}\alpha \right) }\int\limits_{-\infty}^\infty \frac{{\rm sgn}(x-t)}{|t-x|^{1-\alpha}}f(t)dt,\quad {\rm Re}\,\alpha>0,\quad \alpha\neq 2,4,6,...
\end{equation}

For $0<\alpha<1$ operators inverse to $I^\alpha$ and $H^\alpha$ may be constructed in the following forms:
$$
((I^\alpha)^{-1} f)(x)=\frac{1}{2\Gamma(-\alpha)\cos\left(\frac{\pi}{2}\alpha \right) }\int\limits_{-\infty}^\infty \frac{f(x-t)-f(x)}{|t|^{1+\alpha}}dt,
$$
$$
((H^\alpha)^{-1} f)(x)=\frac{1}{2\Gamma(-\alpha)\sin\left(\frac{\pi}{2}\alpha \right) }\int\limits_{-\infty}^\infty \frac{f(x-t)-f(x)}{|t|^{1+\alpha}}{\rm sgn}\,t\,dt.
$$
Next we introduce the convolution 
operator
$$ 
(G^\alpha f)(x)=\int\limits_{-\infty}^\infty G_\alpha(x-t)f(t)dt,
$$
which is defined using the Fourier transform by the equality 
\begin{equation} \label{Bessel02} 
	F[G^\alpha f](x)=\frac{1}{(1+|x|^2)^{\alpha/2}}F[f](x),\qquad {\rm Re}(\alpha)>0.
\end{equation} 
The function $G_\alpha(x)$, whose Fourier transform is $(1+|x|^2)^{-\alpha/2}$, is evaluated in terms of Bessel functions. 
This is why the operator $G_\alpha(x)$ is referred to as a Bessel fractional integration operator or Bessel potential.

Bessel potential or Bessel fractional integral can be written in the form
$$
(G^\alpha f)(x)=\frac{2^{\frac{1-\alpha }{2}}}{\sqrt{\pi } \Gamma \left(\frac{\alpha }{2}\right)} \int\limits_{-\infty }^{\infty } f(t) |x-t| ^{\frac{\alpha -1}{2}} K_{\frac{1-\alpha
	}{2}}(|x-t| ) \, dt,\qquad {\rm Re}(\alpha)>0.
$$

\subsubsection{Fractional integro-differentiation of analytic functions}

Analytic functions are functions that can be represented by convergent power series in a neighborhood of each point in their domain. Fractional integro-differentiation of analytic functions is particularly relevant in the field of fractional calculus, because we can use series expansion. 

Let the function $f(z)=\sum\limits_{k=0}^\infty f_k z^k$, 
be analytic in the unit disc, so 
\begin{equation}\label{GL}
	(\mathfrak{D}^\alpha_0 f)(z)=z^{-\alpha}\sum\limits_{k=0}^\infty  \frac{\Gamma(k+1)}{\Gamma(k-\alpha+1)}f_kz^k.
\end{equation}
A natural way to generalize the \eqref{GL} is to replace the factor $ \frac{\Gamma(k+1)}{\Gamma(k-\alpha+1)}$
in \eqref{GL} by a more general one. One of such generalization is Gel'fond--Leont'ev differentiation of the form (see \cite{Gelfond}) 
\begin{equation}\label{GL01}
	\mathfrak{D}^n (a;f)= \sum\limits_{k=n}^\infty  \frac{a_{k-n}}{a_k}f_kz^{k-n},
\end{equation}
where $a(z)=\sum\limits_{k=0}^\infty a_k z^k$. The operator in \eqref{GL01} is called the Gel'fond--Leont'ev operator of 
generalized differentiation. It is obvious that $\mathfrak{D}^n (a;f)=\frac{d^nf}{dz^n}$  in the case 
$a(z) = e^z$. 

The operator 
$$
\mathfrak{I}^n (a;f)= \sum\limits_{k=0}^\infty  \frac{a_{k+n}}{a_k}f_kz^{k+n},
$$ 
which is the right inverse to \eqref{GL01}, will be called the Gel'fond--Leont'ev  operator
of generalized integration. 
Operators $\mathfrak{D}^n$ and $\mathfrak{I}^n$ are introduced as direct 
generalizations concerning integer order $n$ of integro-differentiation, they contain 
that for fractional order as well. To show this, let us consider the following special 
case when $a(z)=E_\alpha(z)=\sum\limits_{k=0}^{\infty}\frac{z^k}{\Gamma(\alpha k+1)}$, $\alpha>0$
is the Mittag-Leffler function (see Ch. 2, \cite{Kir1} and \cite{Gorenflo}). The corresponding operator for generalized integration of order 
$n=1$ is 
$$
(\euscr{I}_\alpha f)(z)=\mathfrak{I}^1(E_{1/\alpha};f)=\sum\limits_{k=0}^\infty  \frac{\Gamma(\alpha k+1)}{\Gamma(\alpha k+\alpha+1)}f_kz^{k+1}.
$$
This operator admits 
the following integral representation 
$$
(\euscr{I}_\alpha f)(z)=\frac{1}{\Gamma(\alpha)}\int\limits_0^1 (1-t)^{\alpha-1}f(zt^\alpha)dt.
$$

Corresponding differential operator for $0<\alpha<1$ is
$$
(\euscr{D}_\alpha f)(z)=\frac{1}{\Gamma(1-\alpha)}\frac{1}{g'(z)}\frac{d}{dz}\int\limits_0^z \frac{f(t)g'(t)}{[g(z)-g(t)]^\alpha}dt,
$$
where $g(z)=z^{1/\alpha}$. 
The operator $\euscr{D}_\alpha$ corresponds to the expansion 
$$
(\euscr{D}_\alpha f)(z)=\sum\limits_{k=1}^\infty \frac{\Gamma(\alpha k+1)}{\Gamma(\alpha k+1-\alpha)}f_kz^{k-1}.
$$

If we consider the analytic in the unit disc function $b(z)=\sum\limits_{k=0}^{\infty}b_k z^k$ then 
the Hadamard product composition of functions $b(z)$ and 
$f(z)$ is
\begin{equation}\label{HP} 
	\euscr{D}\{b;f\}=b\circ f=\sum\limits_{k=0}^\infty  b_kf_kz^k.
\end{equation}
The operator \eqref{HP} is a very wide generalization of the differentiation. Under the assumption $b_k\rightarrow\infty$
a generalized integration  is
\begin{equation}\label{HP1} 
	\euscr{I}\{b;f\}=\sum\limits_{k=0}^\infty  \frac{f_k}{b_k}\,z^k.
\end{equation}

Choosing various functions $b(z)$ in \eqref{HP} and \eqref{HP1} we obtain integro-differentiation 
operations of various types. Then $b(z)=\frac{\Gamma(\alpha+1)}{(1-z)^{\alpha+1}}$ \eqref{HP} gives Riemann-Liouville fractional differentiation of the function $z^\alpha f(z)$. If we take $b(z)=\frac{\Gamma(\alpha+1)z}{(1-z)^{\alpha+1}}$ in \eqref{HP}
we obtain the Ruscheweyh fractional derivative. 
Function $b(z)=\sum\limits_{k=1}^\infty(ik)^\alpha z^k$ gives the fractional differentiation by Weyl.

\subsubsection{Weyl fractional derivative of a periodic function}

The Weyl fractional derivative of a periodic function is a mathematical concept that extends the traditional notion of differentiation to non-integer orders for periodic functions.

Let $f(x)$ be a $2\pi$-periodic function on $\mathbb{R}$ and let 
$$
f(x)\sim\sum\limits_{k=-\infty}^{\infty} f_ke^{-ikx},\qquad f_k=\frac{1}{2\pi}\int\limits_0^{2\pi}e^{-ikx}f(x)dx
$$
be its Fourier series. 
Here we will consider functions having zero mean value: 
$$
\int\limits_0^{2\pi}f(x)dx=0.
$$

For periodic functions definition  fractional integro-differentiation, 
suggested by Weyl, is   
\begin{equation}	\label{weil01}
	(I_{+}^{(\alpha)}f)(x)\sim \sum\limits_{k=-\infty}^\infty (ik)^{-\alpha}f_ke^{ikx},\quad
	f_k=\frac{1}{2\pi}\int\limits_0^{2\pi}e^{-ikt}f(t)dt,\quad f_0=0.
\end{equation}
Similarly fractional differentiation is defined: 
\begin{equation}	\label{weil02}
	(\mathcal{D}_{+}^{(\alpha)}f)(x)\sim \sum\limits_{k=-\infty}^\infty (ik)^{\alpha}f_ke^{ikx},\quad
	f_k=\frac{1}{2\pi}\int\limits_0^{2\pi}e^{-ikt}f(t)dt,\quad f_0=0.
\end{equation}

The definition \eqref{weil01} may be interpreted as 
\begin{equation}	\label{weil03}
	(I_{+}^{(\alpha)}f)(x)=\frac{1}{2\pi}\int\limits_0^{2\pi}f(x-t)\Psi_{+}^\alpha(t)dt,\qquad \alpha>0,
\end{equation}
where   
\begin{equation}	\label{weil04}
	\Psi_{+}^\alpha(t)=2\sum\limits_{k=1}^\infty\frac{\cos(kt-\alpha\pi/2)}{k^\alpha}= e^{-\frac{1}{2} i \pi  \alpha } \left(e^{i \pi  \alpha } \text{Li}_{\alpha }\left(e^{-i
		t}\right)+\text{Li}_{\alpha }\left(e^{i t}\right)\right),
\end{equation}
$\text{Li}_{\alpha }(z)=\sum\limits_{k=1}^\infty\frac{z^k}{k^\alpha}$ is the fractional order polylogarithm  function.
The right-hand side in \eqref{weil03} 
is called the Weyl fractional integral of order $\alpha$. 

The Marchaud--Weyl derivative is defined as
\begin{equation}	\label{weil06}
	(D_{+}^{(\alpha)}f)(x)=\frac{1}{2\pi}
	\int\limits_0^{2\pi}(f(x)-f(x-t))\frac{d}{dt}\Psi_{+}^\alpha(t)dt,\qquad 0<\alpha<1.
\end{equation}

It is known  
that for $2\pi$-periodic function $f$, such that $f\in L_1(0,2\pi)$ and $\int\limits_0^{2\pi}f(x)dx=0$ the Weyl fractional integrals 
$I_{+}^{(\alpha)}$ for $0<\alpha<1$ 
coincide with the Liouville integrals on the real line $I_{+}^{\alpha}$ \eqref{Eq05}:
$(I_{+}^{(\alpha)} f )(x)=(I_{+}^{\alpha} f )(x)$.
Weyl approach to fractional derivative of a periodic function was considered in \cite{SamkKilbMari1993}.

\subsection{Generalizations of fractional derivatives and applied aspects}\label{GFD}

\subsubsection{Operators with power-logarithmic kernels. Fractional powers of operators}

Operators with power-logarithmic kernels are a class of integral operators that combine power and logarithmic singularities in their kernels. These operators arise in various mathematical contexts, particularly when studying integral equations of the first kind with power-logarithmic kernels. Fractional powers of operators extend the concept of raising an operator to a power.

Let $[a,b]$ be a segment,  $\alpha>0$, $\beta\geq 0$, $\gamma>b-a$.   One of the direct generalizations of the fractional integral $I^\alpha_{a+}$   defined by  \eqref{Eq01} has the form
\begin{equation}	\label{Log01}
	(I^{(\alpha,\beta)}_{a+} f)(x)=\frac{1}{\Gamma(\alpha)}\int\limits_a^{x}\ln^\beta\left( \frac{\gamma}{x-t}\right) {(x-t)^{\alpha-1}}f(t)dt.
\end{equation}
Operators  with a power-logarithmic kernel \eqref{Log01}
were studied in \cite{Kilbas1} for integer $\beta$ and $\gamma=1$. 
In the general case these operators were considered in \cite{Kilbas6}.

Riemann-Liouville fractional integro-differentiation is formally a fractional power 
$\left(\frac{d}{dx} \right)^\alpha$ of the differentiation operator $\frac{d}{dx}$ and is invariant relative to translation 
if considered on the whole axis. 
Hadamard \cite{Hadamard} suggested a construction of fractional integro-differentiation 
which is a fractional power of the type $\left(x\frac{d}{dx} \right)^\alpha$. This construction is well suited 
to the case of the half-axis, and is invariant relative to dilation.

The Hadamard fractional integral has the form   
\begin{equation} \label{Had01} 
	(\mathfrak{F}^{\alpha}_{+}f)(x)=\frac{1}{\Gamma(\alpha)}\int\limits_0^x\frac{f(t)dt}{t\left(\ln\frac{x}{t}\right)^{1-\alpha}},\qquad x>0,\qquad\alpha>0.
\end{equation}

It is easily seen that operator $\mathfrak{F}^{\alpha}_{+}$ is connected with Liouville operators $I^\alpha_{+}$ of the form \eqref{Eq05}  
$$
(\mathfrak{F}^{\alpha}_{+}f)(x)=A^{-1}I^\alpha_{+}Af,\qquad (Af)(x)=f(e^x).
$$

For $0<\alpha<1$ Hadamard fractional derivative  has the form 
$$
(\mathfrak{D}^\alpha_{+}f)(x)=\frac{1}{\Gamma(1-\alpha)}\,x\frac{d}{dx} 
\int\limits_0^x\frac{f(t)dt}{t\left(\ln\frac{x}{t}\right)^{\alpha}}.
$$
We may also consider Hadamard fractional integral and derivative on a finite segment $[a,b]$.

Next we briefly  consider the fractional powers $(B_\gamma)^{\alpha}, \alpha\in \mathbb{R}$ of differential Bessel operator in the form
\begin{equation}\label{Bess}
	B_\gamma= D^2+\frac{\gamma}{x}D=\frac{1}{x^\gamma}\frac{d}{dx}x^\gamma\frac{d}{dx},\qquad \gamma\geq 0, \qquad D:= \frac{d}{dx}.
\end{equation}
For fractional powers of \eqref{Bess}, explicit formulas were derived in \cite{McBArt} as compositions of simpler operators. 
An important step was taken in \cite{Ida}, where explicit definitions were obtained in terms of the Gauss hypergeometric functions, with various applications to PDEs. The most comprehensive study was conducted by I. Dimovski and V. Kiryakova \cite{Dim66,Dim68,DimKir,Kir1} for a more general class of hyper-Bessel differential operators related to the Obrechkoff integral transform.

The  left--sided fractional Bessel integral  $B_{\gamma,a+}^{-\alpha}$  on a segment $[a,b]$ for a function $f{\in}L_1(a,b)$, $a,b\in (0,\infty)$, is defined by  formula
\begin{align}\label{Bess222}
	&(B_{\gamma,a+}^{-\alpha}f)(x)=(IB_{\gamma,a+}^{\alpha}\,f)(x)\nonumber\\
	&=\frac{1}{\Gamma(2\alpha)}\int\limits_a^x\left(\frac{y}{x}\right)^\gamma\left(\frac{x^2{-}y^2}{2x}\right)^{2\alpha-1}\,_2F_1\left(\alpha{+}\frac{\gamma{-}1}{2},\alpha;2\alpha;1{-}\frac{y^2}{x^2}\right)f(y)dy.
\end{align}

Let $\alpha>0$, $n=[\alpha]+1$, $f{\in}L_1(a,b)$, $IB_{\gamma,b-}^{n-\alpha}f,IB_{\gamma,a+}^{n-\alpha}f{\in} C^{2n}(a,b)$ and even. 
The left-sided   fractional Bessel derivatives on a segment of the Riemann-Liouville type for $\alpha\neq 0,1,2,...$ is defined as
$$
(B_{\gamma,a+}^\alpha f)(x)=(DB_{\gamma,a+}^\alpha f)(x)=B_\gamma^n(IB_{\gamma,a+}^{n-\alpha}f)(x),\qquad n=[\alpha]+1.
$$
When
$\alpha=n\in\mathbb{N}\cup\{0\}$, then
$$
(B_{\gamma,a+}^0 f)(x)=f(x),\qquad (B_{\gamma,a+}^n f)(x)=B_{\gamma}^nf(x),
$$
where $B_{\gamma}^n$ is an iterated Bessel operator \eqref{Bess}.
Fractional Bessel integrals and derivatives were studied in \cite{SitSh1}.

\subsubsection{Erd\'{e}lyi-Kober-type operators and fractional integrals and derivatives of a function with respect to another function}

The Erd\'{e}lyi-Kober-type operators   often arise in the context of transmutation operators, integral transforms and functional equations, playing a crucial role in the analysis of differential equations and in the study of special functions such as Bessel functions, hypergeometric functions, and others. Fractional integrals and derivatives with respect to another function extend the concept of fractional calculus to a more generalized setting, where the differentiation or integration operation is performed with respect to a given function rather than a variable. Some basic results concerning Erd\'{e}lyi-Kober operators can be found in the book by Sneddon
\cite{Sneddon}.

Let $0\leq a<x<b\leq \infty$ for any $\sigma\in\mathbb{R}$ or $-\infty\leq a<x<b\leq \infty$ for $\sigma >0$.
Erd\'{e}lyi-Kober-type operators are  
\begin{equation}\label{EK01}
	I_{a+;\,\sigma,\,\eta}^{\alpha} f(x) = \frac{\sigma x^{-\sigma(\alpha+\eta)}}{\Gamma(\alpha)}
	\int\limits_a^x (x^\sigma-t^\sigma)^{\alpha-1}t^{\sigma\eta+\sigma-1}f(t)\,dt,
\end{equation}
for $\alpha>0$, and
\begin{equation}\label{EK02}
	I_{a+;\,\sigma,\,\eta}^{\alpha} f(x) =x^{-\sigma(\alpha+\eta)}\left(\frac{d}{\sigma x^{\sigma-1}dx}\right)^n x^{\sigma(\alpha+n+\eta)}I_{a+;\,\sigma,\,\eta}^{\alpha+n} f(x),
\end{equation}
for $\alpha>-n$, $n\in\mathbb{N}$.

After the change of variables $x^\sigma=y$, $t^\sigma=\tau$ \eqref{EK01}--\eqref{EK02}  are reduced to the 
usual Riemann-Liouville fractional integrals and derivatives 
\begin{equation}\label{EK07}
	I_{a+;\,\sigma,\,\eta}^{\alpha} f(x) = y^{-\alpha-\eta}(I_{a^\sigma+}^{\alpha}\varphi)(y),\qquad \varphi(y)=y^\eta f(x),\qquad x^\sigma=y.
\end{equation}

Erd\'{e}lyi-Kober operators are essential and important in transmutation theory. For example, the most well--known transmutations of Sonine and Poisson are of this class when $\sigma=2$ (see for details \cite{KS1,BookSSh}). 
In the monographs by Sneddon \cite{Sneddon} and by Kiryakova \cite{Kir1} a comprehensive theory of these  operators was given.
Important properties of Erd\'{e}lyi-Kober operators were studied in the  monographs \cite{KS1,BookSShElsevier}.

Let ${\rm Re}\, \alpha>0$, $-\infty\leq a<b\leq\infty$.
The fractional integral of a function $f$ with respect to another function  $g$ is  (see \cite{Holmgren})
\begin{equation}
	\label{163}
	(I_{a+,g}^{\alpha}f)(x)=\frac{1}{\Gamma(\alpha)}\int\limits_a^x \left( g(x)-g(t)\right)^{\alpha-1}g'(t)f(t)d\,t, 
\end{equation}
where $g(x)$ is a strictly increasing function.
Fractional integrals of a function by another in the complex plane were studied by 
Osler \cite{Osler,Osler1}.
Integral \eqref{163} is defined for every function $f(t)\in L_1(a,b)$ and for any monotone function $g(t)$, 
having a continuous derivative. 

If $g'(x)\neq0$, $a\leq x\leq b$, then  operators $I_{a+,g}^{\alpha}$, $I_{b-,g}^{\alpha}$ are expressed via the usual Riemann-Liouville (or Liouville) fractional integration after the corresponding changes of variables
$$
I_{a+,g}^{\alpha}f=QT^\alpha_{c+}Q^{-1}f,\quad I_{b-,g}^{\alpha}f=QT^\alpha_{d-}Q^{-1}f,\quad c=g(a),\quad d=g(b),
$$
where $(Qf)(x)=f[g(x)]$.

The fractional derivative of a function $f$ with respect to another function  $g$ of the order $\alpha\in(0,1)$ is
\begin{equation}
	\label{165}
	(D_{a+,g}^{\alpha}f)(x)=\frac{1}{\Gamma(1-\alpha)}\frac{1}{g'(x)}\frac{d}{dx}\int\limits_a^x \frac{f(t)}{\left(g(x)-g(t)\right)^{\alpha}}\,g'(t)d\,t.
\end{equation}

It is also possible to consider the  Marchaud and other forms of the fractional derivatives of a function $f$ with respect to another function  $g$.

The operational calculus approach for fractional calculus with respect to functions is given in \cite{Fahad}.

Riemann-Liouville fractional integral is obtained by choosing $g(x)=x$ in \eqref{163}. If we take in \eqref{163} the function  $g(x)=x^\sigma$ we obtain
Erd\'{e}lyi-Kober type operator \eqref{EK01}; if we take $g(x)=\ln x$ in \eqref{163} we get Hadamard fractional integral, and a choice $g(x)=\exp(-x)$ with its applications was considered in \cite{DN58}.

As A.M. Djrbashian pointed out to us operators of fractional derivatives of a function with respect to another function \eqref{163} even in some more general setting were introduced and studied by his father M.M. Djrbashian, cf. \cite{DZ1,DZ2,DN3}. In these papers were studied integral representations of this operator class, their inversion and corresponding integro-differential equations of fractional order.

\subsubsection{Variable-order differential operators and sequential fractional derivatives}

In \cite{SR} fractional integration and differentiation
when the order is not a constant but a function was presented. In order to define variable-order differential operators in \cite{SR} two ways were given: the first way is a direct one and the
second uses Fourier transforms. For the Riemann--Liouville fractional integrals of variable order, we have for ${\rm Re}(\alpha(x))>0$, $ x>a$ (see \cite{SR}, formula 2)
\begin{equation}\label{RLVO001}
	(I^{\alpha(x)}_{a+}f)(x)=\frac{1}{\Gamma(\alpha(x))}\int\limits_a^x(x-y)^{\alpha(x)-1}f(y)dy.
\end{equation}
If $a>-\infty$ we have the Riemann definition and for $a=-\infty$ we have the Liouville
definition.

The Riemann-Liouville derivative is also extended to the case of variable order  (see \cite{SR}, formula 3):
\begin{equation}\label{RLVO002}
	(D^{\alpha(x)}_{a+}f)(x)=\frac{1}{\Gamma(1-\alpha(x))}\frac{d}{dx}
	\int\limits_a^x(x-y)^{-\alpha(x)}f(y)dy,
\end{equation}
where $0<{\rm Re}(\alpha(x))<1$, $x>a$.

We observe that the fractional operators \eqref{RLVO001} and \eqref{RLVO002} 
are not inverse to each other
as in the case of constant order, 
as it will be seen below. So, it will not be correct
to introduce $D^{-\alpha(x)}_{a+}$ as $[D^{\alpha(x)}_{a+}]^{-1}$.

For $0\leq{\rm Re}(\alpha(x))<1$ Coimbra derivative \cite{Coimbra0} is
$$
(\,^{CV}D^{\alpha(x)}_0f)(x)=\frac{1}{{\Gamma(1-\alpha(x))}}\int\limits_0^x \frac{f'(y)}{(x-y)^{\alpha (x)}}\ dy+\frac{f(0+)-f(0-)}{x^{\alpha (x)}(\Gamma(1-\alpha(x)))}.
$$

Now let us consider fractional order  depending of two variables $\alpha(x,y)$.
Let $f:[a,b]\rightarrow \mathbb{R}$ be a function.
The left--sided Riemann--Liouville fractional integral of
order $\alpha(\cdot,\cdot)$ is defined by (see \cite{Malinowska})
\begin{equation}\label{RLVO01}
	(\,_aI^{\alpha(\cdot,\cdot)}_xf)(x)=\int\limits_a^x\frac{(x-y)^{\alpha(x,y)-1}}{{\Gamma(\alpha(x,y))}}f(y)dy,\qquad x>a.
\end{equation}

Left--sided Riemann--Liouville derivative of   order $\alpha(x,y)\in(0,1)$ is
\begin{equation}\label{RLVO03}
	(\,_aD^{\alpha(\cdot,\cdot)}_xf)(x)=\frac{d}{dx}\int\limits_a^x\frac{(x-y)^{-\alpha(x,y)}}{{\Gamma(1-\alpha(x,y))}}f(y)dy,\qquad x>a.
\end{equation}
Similarly introduced are Caputo and Marchaud derivatives of variable fractional order (see \cite{SamkoVO,SR}).

Next we consider so called sequential derivatives.
Djrbashian--Nersesyan fractional derivatives, associated with a sequence 
$\{\gamma_0,\,\gamma_1,\,\dots,\,\gamma_m\}$
of order $\sigma$, where
$\sigma=\gamma_0+\gamma_1+...+\gamma_m$,
are defined by 
\begin{equation}
	\label{comp}
	D_{DN}^\sigma=D^{\gamma_0} D^{\gamma_1} \cdots D^{\gamma_m},
\end{equation}
where $D^{\gamma_k}$ are fractional integrals and derivatives of Riemann--Liouville with some endpoint.
These operators were introduced in \cite{DN1,DN2,DN3,DN3E}  and then studied and applied in \cite{DZ7}. The original definitions demand $-1{\leq}\gamma_0{\leq} 0,$ $ 0{\leq}\gamma_k{\leq} 1,$ $1{\leq} k{\leq} m$, as in the above papers integrodifferential equations under such conditions were studied for operators \eqref{comp}. But Djrbashian--Nersesyan fractional operators may be defined and considered for any parameters $\gamma_k$ if appropriate definitions of Riemann--Liouville operators are used.

Djrbashian--Nersesyan operators were patterns for introducing in the book of Miller and Ross \cite{Miller} more general sequential operators of fractional integrodifferentiation for which compositions in definitions of the form \eqref{comp} are consisted of any fractional operators, cf. an useful discussion in \cite{Podlubny}.

\subsubsection{Generalized fractional integrals}

There are a lot of generalizations of fractional integro--differential operators  connected with combinations and compositions of more standard fractional operators. For example, averaged or distributed order fractional operator, associated with any given fractional operator $R^t$,  is introduced by the next formula 
\begin{equation}
	\label{M1}
	I_{MR}^{(a,b)}f=\int\limits_a^b R^t f(t)d\,t,
\end{equation}
where $R^t$ is a given fractional operator of order $t$ of any kind.
In particular, when $R^t$ represents the fractional Riemann–Liouville operator, terms such as "continued"\,
or "distributed-order"\, fractional integrals are often used. Such operators were studied in \cite{Psh1,Psh2}.

Fractional integrals containing the Gaussian hypergeometric function in the kernel are the Saigo fractional integrals (see \cite{Saigo}) which are defined by
\begin{equation}\label{Saigo1}
	J_{x}^{\gamma,\beta,\eta}f(x)=\frac{1}{\Gamma(\gamma)}\int\limits_x^\infty(t-x)^{\gamma-1}t^{-\gamma-\beta}\,_2F_1\left(\gamma+\beta,-\eta;\gamma;1-\frac{x}{t}\right)f(t)dt,
\end{equation}
and
\begin{equation}\label{Saigo2}
	I_{x}^{\gamma,\beta,\eta}f(x)=\frac{x^{-\gamma-\beta}}{\Gamma(\gamma)}\int\limits_0^x(x-t)^{\gamma-1}\,_2F_1\left(\gamma+\beta,-\eta;\gamma;1-\frac{t}{x}\right)f(t)dt,
\end{equation}
where $\gamma>0,\beta,\theta$ are real numbers. 
Similar class of hypergeometric fractional integrals  have been introduced also by Love.  Generalized fractional calculus operators with more general special functions  in the kernels, as Meijer $G$- and Fox $H$-functions have been also studied, see for example \cite{Kalla1} and \cite{Kir1}.

We consider the case of the interval $[0,1]$. Following Hadamard 
\cite{Hadamard} and M.M. Djrbashian \cite{DZ1,DZ2} we introduce the operator 
\begin{equation} \label{Dz01} 
	(L^{(\omega)}f)(x)=-\int\limits_0^1f(xt)\omega'(t)dt,
\end{equation}
where the function $\omega\in C([0,1])$  is supposed to satisfy the following assumptions: 
\begin{enumerate}
	\item $\omega(x)$ is monotone,
	\item $\omega(0)=1$, $\omega(1)=0$, $\omega(x)\neq0$ as $0<x<1$,
	\item  $\omega'(x)\in L_1(0,1)$.
\end{enumerate}

If $\omega(x)=\frac{(1-x)^\alpha}{\Gamma(1+\alpha)}$, then obviously $(L^{(\omega)}f)(x)=x^{-\alpha}(I^\alpha_{0+}f)(x)$, where $I^\alpha_{0+}$ is \eqref{Eq03}. 

M.M. Djrbashian  (also Mkhitar Dzhrbashjan,  M. M. Jerbashian (see \cite{RogosinDubatovskayaD} about his papers))
in \cite{DZ1,DZ2}  considered the operators $L^{(\omega)}$ in a 
more general form 
$$
(L^{(\omega)}f)(x)=-\frac{d}{dx}\left(x\int\limits_0^1 f(xt)dp(t)\right),\qquad p(t)=t\int\limits_t^1 \frac{\omega(x)}{x^2}dx.
$$

\subsubsection{Applied aspects of fractional derivatives and integrals and numerical calculus}\label{Apl}


The reasons for the large number of applications of fractional calculus are as follows.
\begin{enumerate}
	\item The nonlocality of the fractional derivative makes it possible to use it for mathematical modeling of media with memory.
	\item In many models of physics, biology and medicine, differential equations of fractional orders describe the phenomenon under consideration more accurately, since the order of the fractional derivative gives an additional degree of freedom.
	\item 	Fractional equations describe non-Markovian processes, which opens up new horizons in probability theory and statistics.
\end{enumerate}

The use of fractional calculus in theoretical physics and mechanics is described, for example, in the books  \cite{Herrmann,Uchaikin01,Podlubny,Atanackovic}. The mail idea of using fractional derivatives and integrals in physics is based on nonlocal effects which may occur in space and time. Thus, the fractional  derivative by time in the model is interpreted by physicists as the presence of the memory of the described process. A fractional derivative with respect to a spatial variable indicates the presence of movement restrictions.

The paper \cite{KiryakovaAp} introduces  scientists who started to apply fractional calculus to scientific and engineering problems during the nineteenth and twentieth centuries. How to use the Mellin integral transform in fractional calculus was described in \cite{Luchko0,Luchko}. In \cite{Grigoletto} fundamental theorem of fractional calculus is presented.

The use of fractional derivatives in continuum mechanics is due to the emergence of new polymeric materials that have both the property of viscosity and the property of elasticity.  Describing the behavior of solids using viscoelastic models that contain fractional order operators and relate stresses to strains was began in 1948 from papers  of   Rabotnov Yu. N. (see translation \cite{Rabotnov}) and of  Gerasimov A. N. \cite{Gerasimov}.  Let mention here the paper \cite{Nov} where there are many biographical facts and papers of   Gerasimov A. N., including his  paper \cite{Gerasimov} were given.
Models of viscoelasticity with fractional operator were described by Rossikhin Yu. A. and  Shitikova M. V.
in Encyclopedia of Continuum Mechanics \cite{Rossikhin01}.
The mathematical modeling of a stress-strain state of the viscoelastic periodontal membrane using fractional calculus is given in \cite{Bosiakov}. The fractional derivative of strain with time to the ratio between stress and strain is used in this case to describe the intermediate state of the material between elastic and liquid \cite{ShitikovaR,Scott} and \cite{Handbook}, Vol. 7, 139--158 (2019).  
Monograph \cite{Mainardi} provides an overview of application of fractional calculus and waves in linear viscoelastic media, which includes fractional viscoelastic models, dispersion, dissipation and diffusion.  Diethelm K. in \cite{Diethelm} described numerical methods for the fractional differential equations of viscoelasticity. 

In \cite{CaputoSt} it was noticed that in the modelling of the energy dissipation and dispersion in the propagation of elastic waves,  introduction of memory mechanisms in their constitutive equations is needed. It turned out that the most successful memory mechanism used to represent variance and energy dissipation is a fractional derivative. The fractional derivative  generalizing the stress strain relation of anelastic media was introduced in \cite{Caputo2,CaputoSt}.

Models with fractional order operators used in dynamic problems of rigid body mechanics first appeared in \cite{Gemant,Scott,Gerasimov} and described in reviews \cite{ShitikovaR,ShitikovaR01,Rossikhin,ShitikovaMem}. Fractional calculus in the problems of mechanics of solids was given in \cite{Rossikhin02,Rossikhin04,Rossikhin05} and \cite{Handbook}, Vol. 7, 159--192 (2019).

In problems of thermodynamics, solutions for unsteady flows near the boundary of a semi-infinite region are expressed by combinations of fractional derivatives in \cite{Babenko}. In paper
\cite{Meilanov}  thermodynamics was generalized in terms of fractional derivatives, while the results of traditional thermodynamics by Carnot, Clausius and Helmholtz were obtained in the particular case when the exponent of the fractional derivative is equal to one.  In \cite{Povstenko}  was described  how fractional calculus used in thermoelasticity.

Application of fractional modelling in rheology providing by Scott Blair was presented in \cite{Rogosin}.
V.V.  Uchaikin in  \cite{Uchaikin02} presented fractional models in hydromechanics.

Applications of fractional calculus in electromagnetic 
including  fractional multipoles, fractional solutions for Helmholtz equations, and fractional-image methods were studied and briefly reviewed in \cite{Engheta1,Engheta2}. 

In \cite{Petras} fractional-order chaotic systems were studied.
Namely, the model of chaotic system considered in \cite{Petras} as three single differential equations, where the equations contain the  fractional order derivatives. 

Application of the fractional calculus in the control theory can be found in \cite{PodlubnyCont}. It was  shown that 
for fractional-order systems the most suitable way  is to use fractional-order controllers
involving an integrator and differentiator of fractional order.

The works \cite{Zaslavsky1,Zaslavsky2} deal with the fractional generalization of the diffusion equation, fractality of the chaotic dynamics and kinetics, and also includes material on non-ergodic and non-well-mixing Hamiltonian dynamics.

Physical interpretation of a fractional integral can be found in \cite{Rutman}.
Authors of \cite{Kaminsky} gave interpretation of fractional-order operators in fracture mechanics.

Multidisciplinary presentation of fractional calculus for researchers in fields such as electromagnetism, control engineering, and signal processing were given in the book \cite{Ortigueira00}.
In \cite{Bosiakov01} was overviewed  how to use fractional calculus in biomechanics.

Epidemiological SIR model in which classic first order ODE's were replaced by fractional derivatives is gaining popularity due to Covid 19 (see, for example, \cite{SIR}). The presence of fractional derivatives in the model allows you to increase the number of degrees of freedom.

The book \cite{WestBologna}  shows how the fractional calculus can be used to model the statistical behavior of complex systems,  gives general strategies for understanding wave propagation through random media, the nonlinear response of complex materials and the fluctuations of heat transport in heterogeneous materials.

Review \cite{KolokoltsevFr} is fully described advances in the analysis of the fractional and generalized fractional derivatives  from the probabilistic point of view. The theory of fractional Brownian motion and other long-memory processes is describes in \cite{Mishura}, where was  particularly obtained different forms of the Black-Scholes equation for fractional market.
In \cite{BMS1}  Wiener-transformable markets, where the driving process is given by an adapted transformation of a Wiener process including processes with long memory, like fractional Brownian motion and related processes
were studied.
Probabilistic interpretation of a fractional integration is given in \cite{Stanislavsky}.

\vskip 0.5cm

Fractional derivatives and integrals are complex objects, the calculation of which is often technically difficult. In this case, numerical methods are applied. 
Since the fractional derivative is non-local approximation of it is more complicated  than the integer derivatives. 
C. Li with coauthors 
made a great contribution to the development of the numerical calculus of fractional integrals and derivatives and their applications. 
Almost all existing numerical approximations for fractional integrals and derivatives were systematically examined in \cite{Changpin1,Changpin2,Changpin3}.

Difference schemes for solutions to  equation with fractional time derivative was first used in \cite{Shkhanukov}.
In a series of articles \cite{Podlubny0I,Podlubny0II,Podlubny0III} triangular strip matrices  are used
for approximating fractional derivatives and solving fractional differential equations.  

Monte Carlo method is introduced for numerical evaluation of fractional-order derivatives in \cite{PL1,PL2}. A feature of these works is that the points at which the function is calculated are distributed unevenly and their distribution depends on the order of the derivative.

Group of  scientists  as Igor Podlubny, Ivo Petras, Blas  Vinagre  and others deeply advanced in applications of  fractional derivatives in physics including fractional order controller and  distributed-order dynamic systems (see \cite{Podlubny,PodlubnyCont,PodlubnyCont1,Petras,Jiao,Caponetto,PetrasPodlubny}).

V. Kiryakova devoted Section 9 of the article \cite{Kir2} to information on numerical calculations of special functions of fractional calculus.
Recent results on numerical algorithms for evaluation of special functions like  Mittag--Leffler  and Wright can be found in \cite{Gorenflo2,Diethelm1,PodlubnyML} and in \cite{Handbook}, Vol. 1, 269--296 (2019).
Numerical algorithms and results on the  Wright function  and its special
cases can be found in  \cite{LuchkoW,MainardiW}.

\subsubsection{Publications  after 1987}

After 1987 till nowadays, a lot of monographs by authors from different countries had been published, wherein various aspects of fractional calculus, and its applications were considered.

Nowadays, the list of publications wholly or partly devoted to fractional calculus and its applications is endless.
Information about current state of fractional calculus can be found in
\cite{Oliveira,Zhmakin,RogosinDubatovskaya,Debnath,KiryakovaAp,Diethelm01,MKM,Teodoro,Butzer}.

Digest of articles \cite{Handbook} represents a valuable and reliable reference work. The book \cite{Baleanu} provides mathematical tools for working with fractional models and solving fractional differential equations.
The books \cite{Kokilashvili01,Kokilashvili02} provide  insights into multidimensional fractional integral and differential operators and their associated spaces.

We should mention  the great impact of Virginia Kiryakova with her coauthors. She published a book  about generalized fractional derivatives \cite{Kir1} and a lot of papers, including popularizing fractional derivatives \cite{KiryakovaAp,DimKir}. Also she is an editor-in-chief of the most popular journal specialized in fractional calculus: "Fractional Calculus and Applied Analysis" (FCAA; Fract. Calc. Appl. Anal.).

A more detailed list of books on fractional integro-differential calculus is presented in \cite{Handbook}, Vol. 1,  1--22 (2019).

We have discussed various approaches to fractional integro-differentiation. It is essential to apply these methods to fundamental functions such as $z^\lambda$, $e^z$, etc., which are widely used in Taylor and Fourier series.

As demonstrated, there are diverse methods for constructing fractional derivatives and integrals. While the technique of analytic continuation is well-known, it cannot be directly employed to extend analytically the derivative concerning the order of differentiation, as the initial order of the derivative is restricted to natural numbers only. Instead, one can identify an analytic function that aligns with the $n$-th derivative of a function for natural numbers and subsequently extend this function analytically across the entire or almost entire complex plane. Consequently, fractional derivatives and integrals can be constructed in a non-unique manner.

On the other hand, if two distinct integro-differential operators meet the criteria \ref{cr01}, there exists a set of functions where these operators coincide. For instance, when $\alpha>0$ and $f\in L_1(a,b)$, the left-sided Gr\"unwald--Letnikov fractional integral \eqref{GLII} is equivalent to the Riemann--Liouville fractional integral \eqref{Eq01} when $a=0$ for nearly all $x$ (refer to \cite{Podlubny}).

\section{Fractional-order differentiation in the Wolfram language}\label{WS}

\subsection{Riemann-Liouville-Hadamard integro-differentiation in the Wolfram language}\label{Had}

\subsubsection{Definition} 

We describe the Riemann-Liouville integro-differentiation of an arbitrary function to arbitrary symbolic order $\alpha$ 
in the Wolfram Language. Using techniques described below, this operation, hereafter referred to as "fractional differentiation", 
has been published in the Wolfram Function Repository as \texttt{ResourceFunction["FractionalOrderD"]}. 
Defined via an integral transform, fractional differentiation is an analytic function of $\alpha$ which coincides with the usual 
$\alpha$-th order derivative when $\alpha$ is a positive integer and with repeated indefinite integration for negative integer $\alpha$.

We will use notation  $\frac{d^\alpha f(z)}{z^\alpha}$ for Riemann-Liouville integro-diffirentiation  for all $\alpha\in\mathbb{C}$.

{\bf Definition.} By definition of $\frac{d^\alpha f(z)}{dz^\alpha}$ we put
\begin{equation}\label{FrW}
	\frac{d^\alpha f(z)}{dz^\alpha}=\begin{cases}
		f(z), & \alpha =0; \\
		f^{(\alpha )}(z), & \alpha \in \mathbb{Z}\quad \text{and}\quad \alpha >0; \\
		\underbrace{\int\limits_0^zdt...\int\limits_0^{t}dt\int\limits_0^{t}}_{-\alpha}f(t)dt, & \alpha \in \mathbb{Z}\quad \text{and}\quad \alpha <0;\\ 
		\frac{1}{\Gamma (n-\alpha )}\frac{d^n}{dz^n}\int\limits_0^z \frac{f(t) dt}{(z-t)^{\alpha -n+1}}, & n{=}\lfloor \alpha \rfloor +1\quad \text{and}\quad {\rm Re}(\alpha ){>}0; \\
		\frac{1}{\Gamma (-\alpha )}\int\limits_0^z \frac{f(t)dt}{ (z-t)^{\alpha+1}}, & {\rm Re}(\alpha )<0\quad \text{and}\quad \alpha \notin \mathbb{Z} ;\\
		\frac{1}{\Gamma (1-\alpha )}\frac{d}{dz}\int\limits_0^z\frac{ f(t)dt}{(z-t)^{\alpha}}, & {\rm Re}(\alpha )=0\quad \text{and}\quad {\rm Im}(\alpha )\neq 0,
	\end{cases}
\end{equation}
where in the case of divergent integral we use Hadamard finite part. Such integro-differentiation \eqref{FrW} is called {\it Riemann-Liouville-Hadamard derivative}. 

Above we separated cases of symbolic positive integer $n$-th order derivatives from generic result integro-differentiation of fractional order. In particular for $\alpha=-1, -2,...$ we have
$$
\frac{1}{\Gamma (-\alpha )}\int\limits_0^z \frac{f(t)dt}{ (z-t)^{\alpha+1}}=\underbrace{\int\limits_0^zdt...\int\limits_0^{t}dt\int\limits_0^{t}}_{-\alpha}f(t)dt.
$$
So we can combine the third and fifth formulas. Function \texttt{FractionalOrderD} realized regularized  Riemann-Liouville integro-diffirentiation. That means if  any of the integrals in \eqref{FrW} diverges we use  Hadamard
regularization of this  integral (see Subsection \ref{HR}).

Function \texttt{FractionalOrderD} which calculating $\frac{d^\alpha}{dz^\alpha}$ is presented in  Wolfram Function Repository. In order to get this function we should
write 
\begin{center}
	\texttt{ResourceFunction["FractionalOrderD"]}.
\end{center}

\subsubsection{Simplest examples}

Let us consider how the function \texttt{FractionalOrderD}  acts to simple functions.
For example,
\texttt{ FractionalOrderD[$x^2$,$\{x,\alpha\}$]} gives 
$$
\frac{2x^{2-\alpha}}{\Gamma(3-\alpha)}.
$$
and  \texttt{ FractionalOrderD[$\sin[z]$,$\{z,\alpha\}$]} gives 
$$
\left\{ \begin{array}{ll}
	\sin\left(z+\frac{\pi \alpha}{2}\right)  & \mbox{if $\alpha\in\mathbb{Z},\,\alpha \geq 0$};\\
	\frac{2^{\alpha -1}\sqrt{\pi }z^{1-\alpha }}{\Gamma\left(1-\frac{\alpha }{2}\right)\Gamma\left(\frac{3-\alpha }{2}\right) }	\, _1{F}_2\left(1;1-\frac{\alpha }{2},\frac{3-\alpha }{2};-\frac{z^2}{4}\right) & \mbox{in other cases}.\end{array} \right. 
$$
If we put $\alpha=3$ in previous result we obtain $(\sin(z))^{'''}=-\cos(z)$.

When $\alpha=-3$, $f(z)=e^z$ we obtain
\begin{align*}
	\frac{d^{-3} e^z}{dz^{-3}}&=\frac{1}{\Gamma (3)}\int\limits_0^z  (z-t)^{2}e^tdt\\
	&=\int\limits_0^z\left(\int\limits_0^{t_3}\left( \int\limits_0^{t_2}e^{t_1}dt_1\right)dt_2  \right)dt_3\\
	&=e^z-\frac{z^2}{2}-z-1.
\end{align*}
For $\alpha=-1/2$ and for $\alpha=1/2$ we can write
$$
\frac{d^{-1/2} e^z}{dz^{-1/2}}=\frac{1}{\Gamma ({1/2})}\int\limits_0^z  \frac{e^t}{(z-t)^{1/2}}dt=e^z {\rm erf}(\sqrt{z}),
$$
$$
\frac{d^{1/2} e^z}{dz^{1/2}}=
\frac{1}{\Gamma ({1/2})}\frac{d}{dz}\int\limits_0^z  \frac{e^t}{(z-t)^{1/2}}dt= e^z \left(\frac{\Gamma \left(-\frac{1}{2},z\right)}{2\sqrt{\pi }}+1\right),
$$
where ${\rm erf}(z)$ is the integral of the Gaussian distribution, given by
$${\rm erf}(z)=\frac{2}{\sqrt{\pi}}\int\limits_0^z e^{-t^2}dt,$$ $\Gamma \left(\alpha,z\right)$ is the incomplete Gamma function, given by \begin{equation}\label{IncG}
	\Gamma(\alpha,z)=\int\limits_z^\infty t^{\alpha-1}e^{-t}dt.
\end{equation}

Now we can evaluate \texttt{FractionalOrderD} operator for more then  100000 functions, analytical in variable $z$. Let show some examples starting with the "simplest"\, mathematical functions, involving only 
one or two letters $\frac{1}{z}$, $\sqrt{z}$, $z^b$, $a^z$, $e^z$, $z^z$.
Next, we see "named functions" with head-titles like $\log(z)$, $\sin(z)$, $J_\nu(z)$, $J_z(b)$.
At last we see "composed functions" $\sqrt{z^2}$, $\left(z^a\right)^b$, $a^{z^c}$,
$\arcsin(z^3)$. If we apply differentiation or indefinite integration
$$
\left(\frac{1}{z}\right)'=-\frac{1}{z^2},\qquad  (\sqrt{z})'=\frac{1}{2 \sqrt{z}},\qquad
(z^{b})'=bz^{b-1},
$$
$$
(a^z)'=a^z \log (a),\qquad
(e^z)'=e^z,$$
$$ (z^z)'=z^z (\log
(z)+1),\qquad (\log(z))'=\frac{1}{z},\qquad (\sin(z))'=\cos (z),
$$
$$ (J_\nu(z))'= \frac{1}{2} (J_{v-1}(z)-J_{v+1}(z)),\qquad (\sqrt{z^2})'=\frac{z}{\sqrt{z^2}},$$
$$ ((z^{a})^b)'= a b z^{a-1}
\left(z^a\right)^{b-1},\qquad (a^{z^c})'= c \log (a) z^{c-1} a^{z^c},
$$
$$
(\arcsin(z^3))'=\frac{3 z^2}{\sqrt{1-z^6}},\qquad \int \frac{dz}{z}=\log (z)+C,
$$
$$
\int \sqrt{z}dz=\frac{2 z^{3/2}}{3}+C,\qquad \int z^bdz=\frac{z^{b+1}}{b+1}+C,
$$
$$
\int a^zdz=\frac{a^z}{\log
	(a)}+C,\qquad \int e^zdz=e^z+C, 
$$
$$
\int\log(z)dz=-z \log (z)-z+C,\qquad \int\sin(z)=-\cos (z)+C,
$$
$$\int J_\nu(z)dz=\frac{2^{-\nu-1} z^{\nu+1} \Gamma
	\left(\frac{v+1}{2}\right)}{\Gamma\left(\nu+1 \right)\Gamma\left(\frac{\nu+3}{2} \right) } \,
_1 {F}_2\left(\frac{v+1}{2};v+1,\frac{v+3}{2};-\frac{z^2}{4}\right)+C,
$$
$$
\int\sqrt{z^2}dz=\frac{z \sqrt{z^2}}{2}+C,\qquad
\int (z^{a})^bdz=\frac{z \left(z^a\right)^b}{a b+1}+C,
$$
$$
\int a^{z^c} dz=-\frac{z
	\left(\log (a) \left(-z^c\right)\right)^{-1/c} \Gamma \left(\frac{1}{c},-z^c \log
	(a)\right)}{c}+C,
$$
$$
\int \arcsin(z^3)dz=z\arcsin(z^3)-\frac{3}{4} z^4 \,
_2F_1\left(\frac{1}{2},\frac{2}{3};\frac{5}{3};z^6\right)+C.
$$
Here $C$ is arbitrary constant.
We find that not each integral and even derivative can be evaluated. Here there are no results for
$$
(J_z(b))',\qquad \int z^z dz,\qquad \int
J_z(b) \, dz.
$$
Also in many cases we see special functions as result of integration.

Operation \texttt{FractionalOrderD} evaluates repeatable $\alpha$-integer order derivatives and indefinite integrals and extends results for arbitrary complex or real order $\alpha$. For example:
$$
\frac{d^\alpha}{dz^\alpha}\frac{1}{z}=z^{-\alpha } \left(\frac{1}{\Gamma(1-\alpha)}\, _2{F}_2(1,1;2,1-\alpha ;-z)+G_{1,2}^{1,1}\left(z\left|
\begin{array}{c}
	0 \\
	-1,\alpha  \\
\end{array}
\right.\right)\right),
$$
$$
\frac{d^\alpha\sqrt{z}}{dz^\alpha}=\frac{\sqrt{\pi } z^{\frac{1}{2}-\alpha }}{2 \Gamma \left(\frac{3}{2}-\alpha \right)},\qquad \frac{d^\alpha z^b}{dz^\alpha} =\frac{\Gamma (b+1) z^{b-\alpha }}{\Gamma (b-\alpha +1)},
$$
$$
\frac{d^\alpha a^z}{dz^\alpha}=\begin{cases}
	a^z \log ^{\alpha }(a), & \alpha \in \mathbb{Z}, \alpha \geq 0; \\
	a^z \log ^{\alpha }(a) (1-Q(-\alpha ,z \log (a))), & \alpha \in \mathbb{Z}, \alpha <0; \\
	a^z z^{-\alpha } (z \log (a))^{\alpha } (1-Q(-\alpha ,z \log (a))), & \text{in other cases},
\end{cases}
$$
$$
\frac{d^\alpha e^z}{dz^\alpha}=\begin{cases}
	e^z, & \alpha \in \mathbb{Z}, \alpha \geq 0; \\
	e^z (1-Q(-\alpha ,z)), & \text{in other cases},
\end{cases}
$$
$$
\frac{d^\alpha \sin(z)}{dz^\alpha}=\begin{cases}
	\sin \left(\frac{\pi  \alpha }{2}+z\right), & \alpha \in \mathbb{Z}, \alpha \geq 0; \\
	\frac{z^{1-\alpha } \, _1F_2\left(1;1-\frac{\alpha }{2},\frac{3}{2}-\frac{\alpha }{2};-\frac{z^2}{4}\right)}{\Gamma (2-\alpha
		)}, & \text{in other cases},
\end{cases}
$$
$$
\frac{d^\alpha J_\nu(z)}{dz^\alpha}=\frac{\sqrt{\pi } 2^{\alpha -2 \nu} \Gamma (v+1) z^{\nu-\alpha }}{\Gamma\left(\frac{\nu-\alpha
		+1}{2} \right)\Gamma\left(\frac{\nu-\alpha +2}{2} \right)  } \times
$$
$$	
\times	\, _2 {F}_3\left(\frac{\nu+1}{2},\frac{\nu+2}{2};\frac{\nu-\alpha
	+1}{2},\frac{\nu-\alpha +2}{2},\nu+1;-\frac{z^2}{4}\right),
$$
$$
\frac{d^\alpha \sqrt{z^2}}{dz^\alpha}=\frac{\sqrt{z^2} z^{-\alpha }}{\Gamma (2-\alpha )},\qquad \frac{d^\alpha (z^a)^b}{dz^\alpha}=\frac{z^{-\alpha } \left(z^a\right)^b \Gamma (a b+1)}{\Gamma (a b-\alpha +1)},
$$
$$
\frac{d^\alpha \arcsin(z^3)}{dz^\alpha}=\frac{6 z^{3-\alpha }}{\Gamma (4-\alpha )}\times
$$
$$\times \,
_7F_6\left(\frac{1}{2},\frac{1}{2},\frac{2}{3},\frac{5}{6},1,\frac{7}{6},\frac{4}{3};\frac{2}{3}-\frac{\alpha
}{6},\frac{5}{6}-\frac{\alpha }{6},1-\frac{\alpha }{6},\frac{7}{6}-\frac{\alpha }{6},\frac{4}{3}-\frac{\alpha
}{6},\frac{3}{2}-\frac{\alpha }{6};z^6\right),
$$
where $Q(\alpha,z)=\frac{\Gamma(\alpha,z)}{\Gamma(\alpha)}$ is the regularized incomplete Gamma function, $\Gamma(\alpha,z)$ is the incomplete Gamma function \eqref{IncG}. Values $\frac{d^\alpha z^z}{dz^\alpha}$, $\frac{d^\alpha \log(z)}{dz^\alpha}$, $\frac{d^\alpha J_z(b)}{dz^\alpha}$, $\frac{d^\alpha a^{z^c}}{dz^\alpha}$ 
are also calculated but formulas too much complicated.

We would like to compare  Mathematica operations \texttt{ FractionalOrderD}, \texttt{ D}  and  \texttt{Integrate}.
The	direct computation from definition  
in Mathematica gives
$$
D^7\frac{1}{z^2}= \mathtt{D[1/z^2, {z, 7}]}=-\frac{40320}{z^9},
$$
$$
\frac{d^7}{dz^7}\frac{1}{z^2}=\mathtt{ResourceFunction["FractionalOrderD"][1/z^2, {z, 7}]}=-\frac{40320}{z^9},
$$ 
$$
\underbrace{\int\limits_0^zdt...\int\limits_0^{t}dt\int\limits_0^{t}}_{7}\frac{dt}{t^2}= \mathtt{ Nest\left[\int \# dz \&, 1/z^2, 7\right]}=\frac{137 z^5}{7200}-\frac{1}{120} z^5 \log (z),
$$
\begin{align*}
	\frac{d^{-7}}{dz^{-7}}\frac{1}{z^2}&=\mathtt{ResourceFunction["FractionalOrderD"][1/z^2, {z, -7}]}\\
	&=-\frac{1}{120} z^5 \left(\log (z)-\frac{77}{60}\right)=\frac{137 z^5}{7200}-\frac{1}{120} z^5 \log (z).
\end{align*} 
Therefore  \texttt{ FractionalOrderD} generalizes \texttt{ D}  and  \texttt{Integrate}.

We should mention, that differentiation by parameters of functions (for example, $\frac{d^\alpha J_z(b)}{dz^\alpha}$) as rule can not produce well known functions. But we can say the following:
"The first derivative with respect to an upper "parameter" $a_k$, and all derivatives of symbolic integer order m with respect to a "lower" parameter $b_k$ of the generalized hypergeometric function $\,_pF_q(a_1,...,a_p;b_1,...,b_q;z)$, can be expressed in terms of the Kamp\'{e} de F\'{e}riet hypergeometric function of two variables" (see  \cite{blog}).

Let us describe how \texttt{ FractionalOrderD} works. 
For evaluation of $\frac{d^\alpha}{dz^\alpha}$ we are using three approaches.
\begin{enumerate}
	\item Generic formulas  exist for single simple functions; these are converted into pattern matching rules. Here power series, generalized  power series and Hadamard regularization are used.
	\item Convert function to the Meijer $G$-function, using powerful new  operator \texttt{MeijerGReduce} and rule-based \texttt{MeijerGForm} (applies to hypergeometric type functions). Then use formula for the fractional derivative of the Meijer $G$-function.  
	\item Generic formulas for fractional differentiation, which produce Appell function $F_1$ of two variables (in future Lauricella function $D$ of several variables and Humbert function $\Phi$ of two variables will be used).
\end{enumerate}
In the following sections we consider these three approaches.

\subsection{Calculation of fractional derivatives and integrals by series expansion}\label{HR}

\subsubsection{Basic power--logarithmic examples}

Let us consider the first approach to calculating $\frac{d^\alpha}{dz^\alpha}$. 
The \texttt{FractionalOrderD} function allows us to find all these and many other fractional derivatives because
this operation applied to each term of Taylor series expansions of all functions near zero. So if
\begin{equation}\label{H01}
	f(z)=z^{b}\sum\limits_{n=0}^\infty c_n z^{n},
\end{equation}
then
\begin{equation}\label{H02}
	\frac{d^\alpha f(z)}{dz^\alpha}=\sum\limits_{n=0}^\infty c_n \frac{d^\alpha z^{b+n}}{dz^\alpha}.
\end{equation}

Sum representations by formula  \eqref{H01} we meet for functions like 
$$e^z =\sum\limits_{n=0}^{\infty }{\frac {z^{n}}{n!}},\qquad 
J_\nu(z)=\sum\limits_{n=0}^{\infty }{\frac {(-1)^{n}}{n!\Gamma (n+\nu +1)}}{\left({\frac {z}{2}}\right)}^{2n+\nu },$$ but sometimes series expansions include $\log(z)$ function as 
in the logarithmic case of $K_0(z)$:
\begin{equation}\label{H003}
	K_0(z)=-\left(\log\left(\frac{z}{2}\right)+{\gamma} \right)I_0(z)+\sum\limits_{n=1}^\infty  \frac{H_n}{(n!)^2}\left({\frac {z}{2}}\right)^{2n},
\end{equation}
with $n$-th harmonic number $H_n=\sum\limits_{k=1}^n\frac{1}{k}$ and $\gamma$  is Euler--Mascheroni constant
\begin{equation}\label{EM}
	\gamma =\lim_{n\to \infty }\left(-\log n+\sum _{k=1}^{n}{\frac {1}{k}}\right)=-\int\limits_{0}^\infty e^{-x}\log(x)dx.
\end{equation}
It means, that we should consider more general series
\begin{equation}\label{H03}
	f_L(z)=z^b \log^k(z) \sum\limits_{n=0}^{\infty } c_n z^n,
\end{equation}
and  evaluate 
for arbitrary $b$, $\alpha$ and integer $k=0,1,2,...$ the following values
\begin{equation}\label{H05}
	\frac{d^\alpha f_L(z)}{dz^\alpha}=\sum\limits_{n=0}^\infty c_n \frac{d^\alpha }{dz^\alpha}z^{b+n}\log^k(z).
\end{equation}

Let $\lambda=b+n$. In order to calculate \eqref{H02} we should find $\frac{d^\alpha z^\lambda}{dz^\alpha}$ by  formula \eqref{FrW}.
Since $\frac{d^0 z^\lambda}{dz^0}=z^\lambda$ and
for $\alpha \in \mathbb{Z}$ and $\alpha\geq 0$ we have $\alpha=0,1,2,...$ and integer derivative gives
\begin{align}\label{H04}
	\frac{d^\alpha z^\lambda}{dz^\alpha}&= \lambda(\lambda-1)...(\lambda-\alpha+1)z^{\lambda -\alpha }=(\lambda-\alpha+1)_\alpha\, z^{\lambda -\alpha }\nonumber\\
	&=(-1)^{\alpha } (-\lambda )_{\alpha }\,z^{\lambda -\alpha },
\end{align}
where
$(\lambda)_\alpha$ is a Pochhammer symbol.

Therefore we obtain
\begin{equation}\label{FrW01}
	\frac{d^\alpha z^\lambda}{dz^\alpha}=\begin{cases}
		z^\lambda, & \alpha =0; \\
		(\lambda-\alpha+1)_\alpha z^{\lambda-\alpha}, & \alpha \in \mathbb{Z}\quad \text{and}\quad \alpha >0; \\
		\frac{1}{\Gamma (n-\alpha )}\frac{d^n}{dz^n}\int\limits_0^z \frac{ t^\lambda dt}{(z-t)^{\alpha -n+1}}, & n{=}\lfloor{\rm Re}(\alpha )\rfloor{+}1\quad \text{and}\quad {\rm Re}(\alpha ){>}0; \\
		\frac{1}{\Gamma(-\alpha )}\int\limits_0^z \frac{ t^\lambda dt}{ (z-t)^{\alpha+1}}, & {\rm Re}(\alpha )<0;\\
		\frac{1}{\Gamma(1-\alpha )}\frac{d}{dz}\int\limits_0^z\frac{t^\lambda dt}{(z-t)^{\alpha}}, & {\rm Re}(\alpha )=0\quad \text{and}\quad {\rm Im}(\alpha )\neq 0.
	\end{cases}
\end{equation}

In order to calculate \eqref{H05} we should find $\frac{d^\alpha}{dz^\alpha}z^{\lambda}\log^k(z)$ by  formula \eqref{FrW}
\begin{equation}\label{FrW02}
	\frac{d^\alpha}{dz^\alpha}z^{\lambda}\log^k(z){=}\begin{cases}
		z^{\lambda}\log^k(z), & \alpha =0; \\
		(z^{\lambda}\log^k(z))^{(\alpha )}, & \alpha \in \mathbb{Z}\quad \text{and}\quad \alpha >0; \\
		\frac{1}{\Gamma (n-\alpha )}\frac{d^n}{dz^n}\int\limits_0^z \frac{t^{\lambda}\log^k(t) dt}{(z-t)^{\alpha -n+1}}, & n{=}\lfloor \alpha \rfloor{+}1\,\, \text{and}\,\, {\rm Re}(\alpha ){>}0; \\
		\frac{1}{\Gamma (-\alpha )}\int\limits_0^z \frac{t^{\lambda}\log^k(t)dt}{ (z-t)^{\alpha+1}}, & {\rm Re}(\alpha )<0;\\
		\frac{1}{\Gamma (1-\alpha )}\frac{d}{dz}\int\limits_0^z\frac{ t^{\lambda}\log^k(t)dt}{(z-t)^{\alpha}}, & {\rm Re}(\alpha )=0\,\, \text{and}\,\, {\rm Im}(\alpha )\neq 0.
	\end{cases}
\end{equation}
Here for $\alpha \in \mathbb{Z}$ and $\alpha >0$
$$
(z^{\lambda}\log^k(z))^{(\alpha )}=
\underset{j=0}{\overset{\alpha}{\sum }}\binom{\alpha}{j} \frac{\Gamma(\lambda+1)}{\Gamma(\lambda-j+1)} z^{\lambda -j} \frac{d ^{\alpha-j}\log
	^k(z)}{d z^{\alpha-j}}.
$$

If some of integrals 
\begin{align}\label{Int}
	\int\limits_0^z \frac{ t^\lambda dt}{ (z-t)^{\alpha+1}},\qquad \int\limits_0^z \frac{ t^\lambda dt}{(z-t)^{\alpha -n+1}},\nonumber\\ \int\limits_0^z \frac{t^{\lambda}\log^k(t)dt}{ (z-t)^{\alpha+1}},\qquad \int\limits_0^z \frac{t^{\lambda}\log^k(t) dt}{(z-t)^{\alpha -n+1}}
\end{align}
in \eqref{FrW01} or \eqref{FrW02}  diverges we take Hadamard
finite part of this  integral.

\subsubsection{Hadamard finite part} 

The concept of the "finite part"\, of a singular integral introduced by Hadamard
based on  dropping some divergent terms and keeping the finite part. 

Let a function $f=f(x)$ be integrable in an interval $\varepsilon<x<A$ for any
$0<\varepsilon$, $\varepsilon<A<\infty$ and the representation 
\begin{equation}\label{Ad}
	\int\limits_{\varepsilon<x<A}f(x)dx=\sum\limits_{k=1}^N a_k\varepsilon^{-\lambda_k}+h\ln{\frac{1}{\varepsilon}}+J_\varepsilon
\end{equation}
hold valid, where $a_k$, $h$, $\lambda_k$ are some constant positive numbers independent of  $A$.
If   the limit 
$\lim\limits_{\varepsilon\rightarrow 0}J_\varepsilon$ exists, then it is called the Hadamard finite part of the singular integral of the function $f$. The function $f=f(x)$ is said to possess the Hadamard property at the origin. The standard notation for the finite part of the Hadamard singular integral is as follows
\begin{equation}\label{p.f.}
	f.p.\int\limits_{x<A}f(x)dx=\lim\limits_{\varepsilon\rightarrow
		0}J_\varepsilon.
\end{equation}
In the case  $h=0$ in the representation \eqref{Ad}, the function $f=f(x$) is said to possess the non-logarithmic-type
Hadamard property at the origin.

For example, let consider the integral 
$$
\int\limits_0^\infty x^{-1}f(x)dx,
$$
where $f(x)$ is an analytic function on the half infinite interval $[0,\infty)$ such that
$f(0)\neq 0$
and $f(x) = O(x^{-\delta})$, $x\rightarrow\infty$ with $\delta>0$. It is easy to see, that this integral is divergent.
But if we consider the integral 
$$
\int\limits_\varepsilon^\infty x^{-1}f(x)dx,\qquad \varepsilon>0,
$$
then by integrating by part we get
$$
\int\limits_\varepsilon^\infty x^{-1}f(x)dx=-f(\varepsilon)\log(\varepsilon)-\int\limits_\varepsilon^\infty \log(x)f'(x)dx.
$$
Then, the limit 
$$
\lim\limits_{\varepsilon\rightarrow 0}\left(\int\limits_\varepsilon^\infty x^{-1}f(x)dx+f(\varepsilon)\log(\varepsilon) \right) 
$$
exists and is finite. This   limit is the Hadamard finite-part of $\int\limits_0^\infty x^{-1}f(x)dx$ and
$$
f.p.\int\limits_0^\infty x^{-1}f(x)dx=\lim\limits_{\varepsilon\rightarrow 0}\left(\int\limits_\varepsilon^\infty x^{-1}f(x)dx+f(\varepsilon)\log(\varepsilon) \right).
$$
In general, for an analytic function $f(x)$ such that $f(0){\neq}0$
and $f(x){=}O(x^{n-\delta-1})$, $x\rightarrow\infty$ with $\delta>0$ the $f.p.$  of the integral $\int\limits_0^\infty x^{-n}f(x)dx$ is defined by
\begin{align}
	f.p.\int\limits_0^\infty x^{-n}f(x)dx&=\lim\limits_{\varepsilon\rightarrow 0}\left(\int\limits_\varepsilon^\infty x^{-n}f(x)dx\right.
	\nonumber\\
	&\left.-\sum\limits_{k=0}^{n-2}\frac{\varepsilon^{k+1-n}}{k!(n-1-k)!}f^{(k)}(0)+\frac{\log(\varepsilon)}{(n-1)!}f^{(n-1)}(0) \right).
\end{align}

\subsubsection{Hadamard finite part for basic power-logarithmic examples}

Let consider how to take  Hadamard
finite part of
$
\int\limits_0^z \frac{t^\lambda dt}{(z-t)^{\alpha+1}}.
$
For the fourth case in   \eqref{FrW01} we have ${\rm Re}(\alpha)<0$ so we do not have a singular point at $t=z$. 
When ${\rm Re}(\lambda)>-1$ we also do not have a singular point at $t=0$ and can just to calculate integral, so
$$
\int\limits_0^z \frac{t^\lambda dt}{(z-t)^{\alpha+1}}=\frac{\Gamma (-\alpha ) \Gamma (\lambda +1)}{\Gamma (\lambda +1-\alpha)}z^{\lambda -\alpha }\qquad\text{for}\qquad{\rm Re}(\lambda)>-1.
$$
If $-2<{\rm Re}(\lambda)<-1$ we can consider
$$
f.p.\int\limits_0^z \frac{t^\lambda dt}{(z-t)^{\alpha+1}}=\int\limits_0^z \frac{t^\lambda}{(z-t)^{\alpha}}\left( \frac{ 1}{z-t}-\frac{1}{z}\right) dt+\frac{\Gamma (1-\alpha ) \Gamma (\lambda +1) }{\Gamma (\lambda +2-\alpha)}z^{\lambda -\alpha }.
$$
Integral 
$$
\int\limits_0^z \frac{t^\lambda}{(z-t)^{\alpha}}\left( \frac{ 1}{z-t}-\frac{1}{z}\right) dt=
\frac{1}{z}\int\limits_0^z \frac{t^{\lambda+1}dt}{(z-t)^{\alpha+1}}=\frac{\Gamma (-\alpha ) \Gamma (\lambda +2) }{\Gamma (\lambda +2-\alpha)}z^{\lambda -\alpha }
$$ converges for ${\rm Re}(\lambda)>-2$, then 
\begin{align*}
	f.p.\int\limits_0^z \frac{t^\lambda dt}{(z-t)^{\alpha+1}}&=\frac{\Gamma (-\alpha ) \Gamma (\lambda +2) }{\Gamma (\lambda +2-\alpha)}z^{\lambda -\alpha }+\frac{\Gamma (1-\alpha ) \Gamma (\lambda +1) }{\Gamma (\lambda +2-\alpha)}z^{\lambda -\alpha }=\\
	&=\left[\Gamma (-\alpha ) \Gamma (\lambda +2)+\Gamma (1-\alpha ) \Gamma (\lambda +1)\right]\frac{z^{\lambda -\alpha }}{\Gamma (\lambda +2-\alpha)}\\
	&=\left[\Gamma (-\alpha )(\lambda +1) \Gamma (\lambda +1)\right.\\
	&\left.+(-\alpha )\Gamma (-\alpha ) \Gamma (\lambda +1)\right]\frac{z^{\lambda -\alpha }}{\Gamma (\lambda +2-\alpha)}\\
	&=\frac{(\lambda +1-\alpha)\Gamma (-\alpha ) \Gamma (\lambda +1)}{\Gamma (\lambda +2-\alpha)}z^{\lambda -\alpha }\\
	&=
	\frac{\Gamma (-\alpha ) \Gamma (\lambda +1)}{\Gamma (\lambda +1-\alpha)}z^{\lambda -\alpha }.
\end{align*}

For $-n-1<{\rm Re}(\lambda)<-n$, $\lambda \neq -1,-2,...,-n$ we will use the next regularisation  
\begin{align*}
	f.p.\int\limits_0^z \frac{t^\lambda dt}{(z-t)^{\alpha+1}}&=\int\limits_0^z \frac{t^\lambda}{(z-t)^{\alpha}}\left( \frac{ 1}{z-t}-\frac{1}{z}\sum\limits_{k=1}^n\frac{t^{k-1}}{z^{k-1}}\right) dt\\
	&+\sum\limits_{k=1}^n\frac{\Gamma (1-\alpha ) \Gamma (\lambda +k)}{\Gamma (\lambda +k+1-\alpha)}z^{\lambda-\alpha}.
\end{align*}
Let us consider the integral
\begin{align*}
	I&=\int\limits_0^z \frac{t^\lambda}{(z-t)^{\alpha}}\left( \frac{ 1}{z-t}-\frac{1}{z}\sum\limits_{k=1}^n\frac{t^{k-1}}{z^{k-1}}\right) dt\\
	&=\int\limits_0^z \frac{t^\lambda}{(z-t)^{\alpha}}\left( \frac{ 1}{z-t}-\frac{1}{z}-\frac{t}{z^2}-..-\frac{t^{n+1}}{z^n}\right) dt\\
	&=\int\limits_0^z \frac{t^\lambda}{(z-t)^{\alpha}}\left( \frac{ z^n-z^{n-1}(z-t)-z^{n-2}t(z-t)-...-t^{n-1}(z-t)}{z^n(z-t)} \right) dt\\
	&=\int\limits_0^z \frac{t^\lambda}{(z-t)^{\alpha}}\left( \frac{ z^n-z^{n}+z^{n-1}t-z^{n-1}t+z^{n-2}t^2-...-zt^{n-1}+t^{n}}{z^n(z-t)} \right) dt\\
	&=\frac{1}{z^n}\int\limits_0^z \frac{t^{\lambda+n}dt}{(z-t)^{\alpha+1}}
	=\frac{\Gamma (-\alpha ) \Gamma (n+\lambda +1) z^{\lambda -\alpha }}{\Gamma (n-\alpha +\lambda +1)}.
\end{align*}
It converges for ${\rm Re}(\lambda)>-n-1$. Then
\begin{align*}
	f.p.\int\limits_0^z \frac{t^\lambda dt}{(z-t)^{\alpha+1}}&=\frac{\Gamma (-\alpha ) \Gamma (n+\lambda +1) }{\Gamma (n-\alpha +\lambda +1)}z^{\lambda -\alpha }\\
	&+\sum\limits_{k=1}^n\frac{\Gamma (1-\alpha ) \Gamma (\lambda +k)}{\Gamma (\lambda +k+1-\alpha)}z^{\lambda-\alpha}\\
	&=\left( \frac{ \Gamma (n+\lambda +1)}{\Gamma (n-\alpha +\lambda +1)}\right.\\
	&\left.-\alpha\sum\limits_{k=1}^n\frac{\Gamma (\lambda +k)}{\Gamma (\lambda +k+1-\alpha)}\right) \Gamma (-\alpha )z^{\lambda-\alpha}.
\end{align*}
Calculating sum we get 
$$
\sum\limits_{k=1}^n\frac{\Gamma (\lambda +k)}{\Gamma (\lambda +k+1-\alpha)}=\frac{1}{\alpha}\left(\frac{\Gamma (n+\lambda +1)}{\Gamma (n-\alpha +\lambda +1)}-\frac{\Gamma (\lambda +1)}{\Gamma (-\alpha +\lambda +1)} \right). 
$$
Therefore for ${\rm Re}(\lambda)>-n-1$, $\lambda\neq -1,-2,...,-n$
$$
f.p.\int\limits_0^z \frac{t^\lambda dt}{(z-t)^{\alpha+1}}
=
\frac{\Gamma (-\alpha ) \Gamma (\lambda +1)}{\Gamma (\lambda +1-\alpha)}z^{\lambda -\alpha }.
$$
So for all $\alpha$ and $\lambda$ such that $\lambda\neq -1,-2,...,-n$ we have
\begin{equation}\label{power}
	\frac{d^\alpha z^\lambda}{dz^\alpha}=\frac{\Gamma (\lambda +1)}{\Gamma (\lambda +1-\alpha)}z^{\lambda -\alpha }.
\end{equation}

Now let consider the case when $\lambda \in \mathbb{Z}$ and $\lambda <0$ i.e. $\lambda= -1,-2,...,-n$. Here we have two variants. The first one it is the case when $\alpha \in \mathbb{Z}$ and $\lambda <\alpha$.
When $\alpha \in \mathbb{Z}$ and  $\alpha<0$ we have $\alpha=-1,-2,...$ and 
$$
\frac{d^{-1} z^\lambda}{dz^{-1}}=f.p.\int\limits_0^z t^{\lambda}dt.
$$
Since for $\varepsilon>0$
$$
\int\limits_\varepsilon^z t^{\lambda}dt=\frac{z^{\lambda+1}}{\lambda+1}-\frac{\varepsilon^{\lambda+1}}{\lambda+1}
$$
then 
$$
\frac{d^{-1} z^\lambda}{dz^{-1}}=f.p.\int\limits_0^z t^{\lambda}dt=\frac{z^{\lambda+1}}{\lambda+1}\qquad \text{for}\qquad \lambda=-2,-3,...
$$
Similarly, for $\alpha=-2,-3,...$ and $\lambda <\alpha$, $\lambda\in\mathbb{Z}$
\begin{align*}
	\frac{d^{\alpha} z^\lambda}{dz^{\alpha}}=f.p.\underbrace{\int\limits_0^zdt...\int\limits_0^{t}dt\int\limits_0^{t}}_{-\alpha}t^\lambda dt\\
	&=\frac{z^{\lambda-\alpha}}{(\lambda+1)(\lambda+2)...(\lambda-\alpha)}.
\end{align*}
Let $\alpha=-n$, $\lambda=-m$. Since $m>n$ we can write 
\begin{align*}
	\frac{1}{(\lambda+1)(\lambda+2)...(\lambda-\alpha)}&=\frac{1}{(\lambda+1)(\lambda+2)...(\lambda+n)}
	\\
	&=\frac{(-1)^n (-1-n-\lambda)!}{(-1-\lambda)!}\\
	&=\frac{\lambda!}{(\lambda-\alpha)!}=(-1)^{\alpha } (-\lambda )_{\alpha }.
\end{align*}
So for all $\alpha$ and $\lambda$ such that $\alpha \in \mathbb{Z}$, $\lambda \in \mathbb{Z}$, $\lambda <0$ and $\lambda <\alpha$   taking into account \eqref{H04} we have
\begin{equation}\label{power01}
	\frac{d^\alpha z^\lambda}{dz^\alpha}=(-1)^{\alpha } (-\lambda )_{\alpha }z^{\lambda-\alpha}.
\end{equation}

In the case $\lambda \in \mathbb{Z}$, $\lambda <0$ and  $\lambda\geq\alpha$ 
the integration of $z^\lambda$ can produce $\log(z)$  	(for negative integer $\lambda$ when $\lambda\geq\alpha$).
For example, when $\lambda=-1$ and $\alpha=-1$
$$
\int\limits_\varepsilon^z \frac{dt}{t}=\log(z)-\log(\varepsilon)
$$  
so  by \eqref{Ad} and \eqref{p.f.} we get
$$
\frac{d^{-1} z^{-1}}{dz^{-1}}= f.p.\int\limits_0^z \frac{dt}{t}=\log(z).
$$
The same regularisation we can obtain  
if we consider  the Laurent series expansion by $\lambda$ at $\lambda=-1$: 
\begin{align*}
	\frac{d^{-1} z^{\lambda}}{dz^{-1}}&=\int\limits_0^z  t^\lambda dt 
	\\
	&=\frac{z^{\lambda +1 }}{\lambda+1}=\frac{1}{\lambda +1}\left[1+(\lambda +1) \log (z)+O\left((\lambda +1)^2\right)\right]\\
	&=\frac{1}{\lambda +1}+\log (z)+O\left(\lambda +1\right).
\end{align*}
Then throwing away the main part $\frac{1}{\lambda +1}$ we get
$$
\frac{d^{-1} z^{-1}}{dz^{-1}}= f.p.\int\limits_0^z \frac{dt}{t}=\lim\limits_{\lambda\rightarrow-1}[\log (z)+O\left(\lambda +1\right)]=\log(z).
$$
In general case when $\lambda=-n$, $n\in\mathbb{Z}$, $\alpha\leq -n$ the Laurent series expansion of \eqref{power} by $\lambda$ at $\lambda=-n$ is 
\begin{align*}
	\frac{\Gamma (\lambda +1)}{\Gamma (\lambda +1-\alpha)}z^{\lambda -\alpha }&=\frac{(-1)^{\lambda-1}z^{\lambda-\alpha}}{(-1-\lambda)!\,\Gamma(1+\lambda-\alpha)(\lambda+n)}\\
	&+\frac{(-1)^{\lambda-1} (\psi(-\lambda)-\psi(\lambda-\alpha+1)+\log(z))}{(-\lambda-1)!\Gamma(\lambda-\alpha+1)} z^{\lambda-\alpha}\\
	&+O\left(\lambda +n\right).
\end{align*}
Here $\psi$ is the diGamma function, given by $\psi(z)=\frac{\Gamma'(z)}{\Gamma(z)}$. 
Then throwing away the main part of the Laurent series expansion  we obtain
\begin{align}
	\frac{d^{\alpha} z^{-n}}{dz^{\alpha}}&= f.p.\int\limits_0^z \frac{t^\lambda dt}{(z-t)^{\alpha+1}}\nonumber\\
	& =\lim\limits_{\lambda\rightarrow-n}\left[\frac{(-1)^{\lambda-1} (\psi(-\lambda)-
		\psi(\lambda-\alpha+1)+\log(z))}{(-\lambda-1)!\Gamma(\lambda-\alpha+1)} z^{\lambda-\alpha}+O\left(\lambda +n\right)\right]\nonumber\\
	&=\frac{(-1)^{-n-1} (\psi(n)-\psi(1-n-\alpha)+\log(z))}{(n-1)!\Gamma(1-n-\alpha)} z^{-n-\alpha}.
\end{align}

So 	we  put
\begin{equation}\label{power02}
	\frac{d^\alpha z^\lambda}{dz^\alpha}=\frac{(-1)^{\lambda-1} (\psi(-\lambda)-\psi(\lambda-\alpha+1)+\log(z))}{(-\lambda-1)!\Gamma(\lambda-\alpha+1)} z^{\lambda-\alpha}
\end{equation}
under certain restrictions.

Finally, combining formulas \eqref{power}, \eqref{power01} and \eqref{power02} we get
\begin{align}\label{powerG}
	\frac{d^\alpha z^\lambda}{dz^\alpha}&\nonumber\\
	&=z^{\lambda -\alpha } \left\{
	\begin{array}{cc}
		& 
		\begin{array}{cc}
			(-1)^{\alpha } (-\lambda )_{\alpha }, & \alpha \in \mathbb{Z}, \lambda \in \mathbb{Z}, \lambda <0, \lambda <\alpha;  \\
			\frac{(-1)^{\lambda -1} (\log (z)+\psi(-\lambda )-\psi(1-\alpha +\lambda))}{(-\lambda -1)! \Gamma (1-\alpha +\lambda)}, & \lambda \in \mathbb{Z}, \lambda <0 ;  \\
			\frac{\Gamma (\lambda +1)}{\Gamma (\lambda+1-\alpha)}, & \text{in other cases}. \\
		\end{array}
		\\
	\end{array}
	\right.
\end{align}
The second integral $\int\limits_0^z \frac{ t^\lambda dt}{(z-t)^{\alpha -n+1}}$ in \eqref{Int} is calculated in same way since ${\rm Re}(\alpha -n)<0$.

We also need to calculate other integrals in \eqref{Int}.
Calculation and taking the Hadamard finite part of  $ \int\limits_0^z \frac{t^{\lambda}\log^k(t)dt}{ (z-t)^{\alpha+1}}$ and $\int\limits_0^z \frac{t^{\lambda}\log^k(t) dt}{(z-t)^{\alpha -n+1}}
$ based on the same ideas as calculation of $\int\limits_0^z \frac{ t^\lambda dt}{(z-t)^{\alpha+1}}$ but much more complicated.  For example, when ${\rm Re}(\alpha )<0$, ${\rm Re}(z)>0$ and ${\rm Re}(\lambda )>-1$ we have
\begin{equation}\label{DLog}
	\frac{d^\alpha }{dz^\alpha}z^\lambda\log(z)=z^{\lambda -\alpha }\frac{\Gamma (\lambda +1)  \left(H_{\lambda }-H_{\lambda -\alpha }+\log
		(z)\right)}{\Gamma (\lambda-\alpha+1)},
\end{equation}
where $H_s= \gamma+\psi(s+1)$ is the harmonic number, $H_n=\sum\limits_{k=1}^n \frac{1}{k}$ for integer $n$,  $\psi(z)=\frac{\Gamma'(z)}{\Gamma(z)}$ is the diGamma function, $\gamma$ is Euler--Mascheroni constant.
However, in general
$$
\frac{d^\alpha }{dz^\alpha}\left( z^\lambda\log(z)\right) =
$$
\begin{equation}\label{DLogG}
	=z^{\lambda -\alpha }\begin{cases}
		(-1)^{\alpha } (-\lambda)_{\alpha }(\psi(-\lambda)-\psi(\alpha -\lambda)+\log (z)); & \alpha,\lambda \in \mathbb{Z},\lambda <0,\lambda <\alpha,  \\
		\begin{split}	\frac{(-1)^{\lambda+1}}{{2\Gamma(-\lambda)\Gamma(\lambda-\alpha+1)}} \biggl(\log^2(z)+\\
			+\frac{\pi^2}{3}-\psi ^{(1)}(\lambda-\alpha+1)-\psi^{(1)}(-\lambda)+\\
			+(\psi(-\lambda)-\psi(\lambda-\alpha+1))\times\\
			\times\left(\psi(-\lambda)-H_{\lambda-\alpha }+
			2\log(z)+\lambda \right)\biggr);
		\end{split} & \lambda\in \mathbb{Z},\lambda <0, \\
		\frac{\Gamma (\lambda+1)\left(H_{\lambda}-H_{\lambda-\alpha }+\log (z)\right)}{\Gamma (\lambda-\alpha+1)}; & \text{in other cases,}
	\end{cases}
\end{equation}
Here 
$\psi ^{(m)}(z)={\frac {\mathrm {d} ^{m+1}}{\mathrm {d} z^{m+1}}}\ln \Gamma (z)
$ is the polyGamma function of order $m$.

Expanding basic functions in Mathematica to sums of the forms \eqref{H01}
or \eqref{H03} allowed us to build fractional derivatives of order $\alpha$ for more than  100000
test cases, involving basic functions and their compositions.

\subsection{The Meijer $G$-function and fractional calculus}\label{MG}

\subsubsection{Definitions and main properties}

All known special functions can be divided on several large groups: 
\begin{itemize}
	\item basic special functions, 
	\item hypergeometric type functions, 
	\item Riemann Zeta and related functions, 
	\item elliptic and related functions, 
	\item number theoretic functions, 
	\item generalized functions and other standard special functions.
\end{itemize}
Basic special functions refer to a class of mathematical functions that have well-established properties and play a fundamental role in various areas of mathematics, physics, and engineering. These functions often arise in the solutions of differential equations, integrals, and other mathematical problems. Some of the most commonly encountered basic special functions include Gamma function, Beta function, Bessel functions,  hypergeometric functions.
Informations about some special functions can be found in \cite{NIST}.

The majority of the most important special functions are the analytic functions, which by definitions typically can be presented through infinite series:
\begin{equation}\label{SG}
	\left(z-z_0\right){}^{\alpha } \log ^k\left(z-z_0\right) \sum _{n=0}^{\infty } c_n\left(z-z_0\right)^n,
\end{equation}
where $k\in\mathbb{Z}$, $k\geq 0$, $c_j,\alpha \in \mathbb{C}$, $j=0,1,2...$.
Such series can be analytically continued (re-expanded) or represented through integrals. 
Very often they are solutions of corresponding differential equations, which came from applications and define special functions.

Hypergeometric type functions can be defined as the functions, which generically can be represented through linear combinations of Meijer $G$-function  which is a very general special function of the form \cite{site4}, see also handbook \cite{Prudnikov}:
\begin{align}\label{G}
	G_{p,q}^{m,n}&=	G_{p,q}^{m,n}\left(z \left| 
	\begin{array}{c}
		a_1,...,a_n,a_{n+1},...,a_p \\
		b_1,...,b_m,b_{m+1},...,b_q \\
	\end{array}\right.
	\right)\nonumber\\
	&={\frac {1}{2\pi i}}\int\limits_{\mathcal{L}}{\frac{\prod\limits_{k=1}^{m}\Gamma (b_{k}-s)\prod\limits_{k=1}^{n}\Gamma (1-a_{k}+s)}{\prod\limits_{k=n+1}^{p}\Gamma (a_{k}-s)\prod\limits_{k=m+1}^{q}\Gamma (1-b_{k}+s)}}\,z^{s}\,ds,
\end{align}
where contour $\mathcal{L}$  represents the path to be followed while integrating. 
This integral is of the so-called Mellin--Barnes type, and may be viewed as a generalization of  inverse Mellin transform.
Information about using of Mellin-Barnes integrals in asymptotic analysis can be found in \cite{Paris}. 
There are three different paths $\mathcal{L}$ of integration:
\begin{itemize}
	\item $\mathcal{L}$ runs from $-i\infty$ to $+i\infty$ 
	so that all poles of $\Gamma(b_i-s)$, $i=1,...,m$, are to the right, 
	and all the poles of  $\Gamma(1-a_k+s)$, $k=1,...,n$, 
	to the left, of $\mathcal{L}$. 
	This contour can be a vertical straight  line $(\gamma-i\infty,\gamma+i\infty)$
	if  ${\rm{Re}}(b_i-a_k)>-1$   for $i =1,...,m$ and
	$k=1,...,n$, (then ${\rm{Re}}(a_k)-1<\gamma<{\rm{Re}}(b_i)$). 
	The integral converges if $p+q<2(m+n)$ and
	$|{\rm arg}\,z|<\left(m+n-\frac{p+q}{2}\right)\pi$.
	If $m+n-\frac{p+q}{2}=0$, then $z$ must be real and positive and additional condition 
	$(q-p)\gamma+\rm{Re}(\mu)<0$,
	$\mu=\sum\limits_{i=1}^q b_i-\sum\limits_{k=1}^p a_k+\frac{p-q}{2}+1$, should be added.  
	\item $\mathcal{L}$ is a loop starting and ending at $+\infty$
	and encircling all poles of
	$\Gamma(b_i-s)$, $i=1,...,m$, once in the negative direction, but none of
	the poles of $\Gamma(1-a_k+s)$, $k=1,...,n$. 
	The integral converges
	if $q\geq 1$ and either $p<q$ or $p=q$ and $|z|<1$ or $q=p$ and $|z|=1$ and $m+n-\frac{p+q}{2}\geq 0$ and $\rm{Re}(\mu)<0$.
	\item $\mathcal{L}$ is a loop starting and ending at $-\infty$
	and encircling all poles of
	$\Gamma(1-a_k+s)$, $k=1,...,n$, once in the positive direction, but none
	of the poles of $\Gamma(b_i-s)$, $i=1,...,m$. 
	The integral converges if
	$p\geq 1$ and either $p>q$ or $p=q$ and $|z|>1$ or $q=p$ and $|z|=1$ and $m+n-\frac{p+q}{2}\geq 0$ and $\rm{Re}(\mu)<0$.
\end{itemize}

Above definition of the Meijer $G$-function holds under the following assumptions:
\begin{itemize}
	\item $0\leq m \leq q$ and $0\leq n\leq p$, where $m$, $n$, $p$ and $q$ are integer numbers,
	\item $a_k-b_j\neq 1,2,3,...$ for $k=1,2,...,n$ and $j=1,2,...,m$, which implies that no pole of any $\Gamma(b_j-s)$,
	$j=1,2,...,m$, coincides with any pole of any $\Gamma(1-a_k+s)$, $k=1,2,...,n$.
	\item $z\neq0$.
\end{itemize}

\begin{rem} 
	A different from \eqref{G} but equivalent to it form was used in the books \cite{Marichev,Prudnikov,MarichevMellin}: 
	\begin{align}\label{G11}
		G_{p,q}^{m,n}&=
		G_{p,q}^{m,n}\left(z \left| 
		\begin{array}{c}
			a_1,...,a_n,a_{n+1},...,a_p \\
			b_1,...,b_m,b_{m+1},...,b_q \\
		\end{array}\right.
		\right)\nonumber\\
		&={\frac {1}{2\pi i}}\int\limits_{\mathcal{L}}{\frac{\prod\limits_{k=1}^{m}\Gamma (b_{k}+s)\prod\limits_{k=1}^{n}\Gamma (1-a_{k}-s)}{\prod\limits_{k=n+1}^{p}\Gamma (a_{k}+s)\prod\limits_{k=m+1}^{q}\Gamma (1-b_{k}-s)}}\,z^{-s}\,ds
	\end{align}
	where
	$$
	m,n,p,q\in \mathbb{N}\cup\{0\},\qquad m\leq q,\qquad n\leq p.
	$$
	Both forms  can be transformed each to other by changing variable of integration $s\rightarrow -s$.
\end{rem}

In the system Mathematica for standard form of the Meijer $G$-function the following notation is used
$$
\mathtt{MeijerG[\{\{a_1,...,a_n\}, \{a_{n+1},...,a_p\}\}],\{\{b_1,...,b_m\}, \{b_{m+1},...,b_q\}\},z]}=
$$
$$
=\mathtt{\frac{1}{2\pi i} {ContourIntegrate}\left(\frac{  \prod\limits_{k=1}^n \Gamma \left(1-a_k-s\right)  \prod\limits_{k=1}^m \Gamma
		\left(b_k+s\right)}{\prod\limits_{k=n+1}^p \Gamma \left(a_k+s\right) \prod\limits_{k=m+1}^q \Gamma
		\left(1-b_k-s\right)}z^{-s},\{s,\mathcal{L}\}\right)};
$$
$$	
\mathtt{ m\in \mathbb{Z}\land m\geq 0\land n\in \mathbb{Z}\land n\geq 0\land p\in
	\mathbb{Z}\land p\geq 0\land q\in \mathbb{Z}\land q\geq 0\land m\leq q\land n\leq p.}
$$

In many classical special functions (like Bessel functions) instead of $z$ the construction $c z^{1/r}$ can be used and under integrand factor 
$z^{-s}$ becomes $(c z^{1/r})^{-s}$ which is not equal to the expression $c^{-s} z^{-s/r}$. 
Correct formulas for these transformation one can find at  Wolfram Function Site, for example in \cite{site}. 

Described properties are supported in Mathematica by command PowerExpand with or without assumption option:

\{PowerExpand[($(z^a)^b$), Assumptions $\rightarrow$ True], 
PowerExpand[$(z^a)^b$]\}

$$
\left\{ e^{2 i \pi  b \left\lfloor \frac{1}{2}-\frac{\rm{Im}(a \log (z))}{2 \pi }\right\rfloor }z^{a b},z^{a b}\right\}
$$
This situation stimulated us to define generalized form of the Meijer $G$-function with additional real parameter $r$, which allowed effectively without restrictions work in the full complex plane with functions like $J_a(z)$.

\begin{rem}
	The generalized form of the Meijer $G$-function with additional real parameter $r$ is defined in \cite{FSite} by similar to \eqref{G} integral
	$$
	G_{p,q}^{m,n}\left(z,r \left| 
	\begin{array}{c}
		a_1,...,a_n,a_{n+1},...,a_p \\
		b_1,...,b_m,b_{m+1},...,b_q \\
	\end{array}\right.
	\right)
	$$
	\begin{equation}\label{Gr}
		={\frac {r}{2\pi i}}\int\limits_{\mathcal{L}}{\frac{\prod\limits_{k=1}^{m}\Gamma (b_{k}+s)\prod\limits_{k=1}^{n}\Gamma (1-a_{k}-s)}{\prod\limits_{k=n+1}^{p}\Gamma (a_{k}+s)\prod\limits_{k=m+1}^{q}\Gamma (1-b_{k}-s)}}\,z^{-\frac{s}{r}}\,ds,
	\end{equation}
	where $r\in \mathbb{R}$, $r\neq 0$, $m\in \mathbb{Z}$, $m\geq 0$, $n\in \mathbb{Z}$, $n\geq 0$, $p\in \mathbb{Z}$, $p\geq 0$,
	$q\in\mathbb{Z}$, $q\geq 0$, $m\leq q$, $n\leq p$.
\end{rem}

Evidently that for default case $r = 1$ we have equality  
$$
G_{p,q}^{m,n}\left(z,1 \left| 
\begin{array}{c}
	a_1,...,a_n,a_{n+1},...,a_p \\
	b_1,...,b_m,b_{m+1},...,b_q \\
\end{array}\right.
\right)=G_{p,q}^{m,n}\left(z \left| 
\begin{array}{c}
	a_1,...,a_n,a_{n+1},...,a_p \\
	b_1,...,b_m,b_{m+1},...,b_q \\
\end{array}\right.
\right).
$$

This Meijer $G$-function with parameter $r$ satisfies two important properties
\cite{site1,site2}:
$$
G_{q,p}^{n,m}\left(\frac{1}{z},r \left| 
\begin{array}{c}
	1-b_1,...,1-b_m,1-b_{m+1},...,1-b_q \\
	1-a_1,...,1-a_n,1-a_{n+1},...,1-a_p \\
\end{array}\right.
\right)
$$
\begin{equation}\label{G111}		
	=G_{p,q}^{m,n}\left(z,r \left| 
	\begin{array}{c}
		a_1,...,a_n,a_{n+1},...,a_p \\
		b_1,...,b_m,b_{m+1},...,b_q \\
	\end{array}\right.
	\right),\,  {z\notin \mathbb{R}},
\end{equation}
$$
G_{p,q}^{m,n}\left(z,r \left| 
\begin{array}{c}
	\alpha+a_1,...,\alpha+a_n,\alpha+a_{n+1},...,\alpha+a_p \\
	\alpha+b_1,...,\alpha+b_m,\alpha+b_{m+1},...,\alpha+b_q \\
\end{array}\right.
\right)
$$
\begin{equation}\label{G1111}	
	=z^{\alpha/r}G_{p,q}^{m,n}\left(z \left| 
	\begin{array}{c}
		a_1,...,a_n,a_{n+1},...,a_p \\
		b_1,...,b_m,b_{m+1},...,b_q \\
	\end{array}\right.
	\right),
\end{equation}

\subsubsection{Representability through Meijer $G$-functions}

Using classical and generalized Meijer $G$-functions we can represent Bessel function $J_\nu(z)$ by formulas
\begin{equation}\label{Bess02}
	J_{\nu }({z})=G_{0,2}^{1,0}\left(\frac{z^2}{4} \left| 
	\begin{array}{c}
		- \\
		\frac{\nu}{2},-\frac{\nu}{2} \\
	\end{array}\right.
	\right),\qquad -\frac{\pi }{2}<\arg (z)\leq \frac{\pi }{2},
\end{equation}
\begin{equation}\label{Bess01}
	J_{\nu }(z)=G_{0,2}^{1,0}\left(\frac{z}{2},\frac{1}{2} \left| 
	\begin{array}{c}
		- \\
		\frac{\nu}{2},-\frac{\nu}{2} \\
	\end{array}\right.
	\right).
\end{equation}
First formula includes parameter $r=1/2$ 
but it is valid for all complex $z$-plane. 
Second formula including classical $G$-function is correct only for half plane.  
For finding last representation through Meijer $G$-function we can use commands  \texttt{ResourceFunction["MeijerGForm"]} and \texttt{MeijerGReduce}:
\begin{center}
	\texttt{{ResourceFunction["MeijerGForm"][BesselJ[$\nu$, z], z], 
			MeijerGReduce[BesselJ[$\nu$, z], z]}}.
\end{center}

For best understanding Meijer $G$-function we can suggest to use several internal commands \texttt{MeijerGInfo}, \texttt{MeijerGToSums}, \texttt{SlaterForm}, which operate with context \texttt{System`MeijerGDump`}, for example,

\texttt{System`MeijerGDump`SlaterForm[MeijerG[\{\{a\}, \{b\}\}, \{\{c\}, \{d\}\}, z, 2/3],
	s]}\\
gives
$$
\frac{2 z^{-s} \Gamma \left(1-a-\frac{2 s}{3}\right) \Gamma \left(c+\frac{2 s}{3}\right)}{3 \Gamma \left(b+\frac{2
		s}{3}\right) \Gamma \left(1-d-\frac{2 s}{3}\right)}.
$$
\texttt{SlaterForm} applied to the Meijer $G$-function reveals the internal structure of the Meijer $G$-function, including the ratio of products of groups with the Gamma function corresponding to their lengths $m$, $n$, $p-n$, $q-m$.

Through Meijer $G$-function we can represent more than 100 named well known functions such as $\log$, $\exp$, $\arctan$, Bessel, Airy and Legendre functions. 
In Mathematica \texttt{TraditionalForm} of these functions form the following list with 127 functions to which we can apply   \texttt{ResourceFunction["MeijerGForm"]}. See below List  \texttt{MeijerG127} in  \texttt{TraditionalForm}
\begin{center}
	\texttt{MeijerG127}=\{$a^z$, $e^z$, $\sqrt{z}$, $z^b$, ${\rm Ai}\,(z)$, ${\rm Ai}'(z)$, ${\rm Bi}(z)$, ${\rm Bi}'(z)$,
	$\pmb{J}_a(z)$,  $\pmb{J}_a^b(z)$,
	${\rm arccos}\,(z)$, ${\rm arcosh}\,(z)$,  ${\rm arccot}\,(z)$, ${\rm arcoth}\,(z)$,
	${\rm arccsc}\,(z)$, ${\rm arcsch}\,(z)$,
	${\rm arcsec}\,(z)$, ${\rm arsech}\,(z)$, 
	${\rm \arcsin}\,(z)$, ${\rm arsinh}\,(z)$,
	${\rm arctan}\,(z)$,  ${\rm arctan}\,(a,z)$,  ${\rm arctan}\,(z,a)$,  ${\rm artanh}\,(z)$,
	$I_a(z)$, $J_a(z)$,  $K_a(z)$, $Y_a(z)$,  $B_z(a,b)$, $B_{(c,z)}(a,b)$, $B_{(z,c)}(a,b)$, $I_z(a,b)$, 
	$I_{(c,z)}(a,b)$, $I_{(z,c)}(a,b)$,
	$R_C(x,z)$, $R_C(z,y)$, $R_E(x,z)$, $R_E(y,z)$, $R_K(x,z)$, 
	$R_K(y,z)$, $T_a(z)$, $U_a(z)$, 
	$\cos (z)$, $\cosh (z)$, $\text{Chi}(z)$, $\text{Ci}(z) $,
	$ F(z)$, $E(z)$, $ K(z)$, $ \text{erf}(z)$, $\text{erf}(a,z)$, $\text{erf}(z,b)$,  $\text{erfc}(z)$, 
	$\text{erfi}(z)$,  $E_a(z)$, $\text{Ei}(z)$, $F_z$, $F_a(z)$, 
	$C(z)$, $F(z)$, $G(z)$,  $S(z)$, $ \Gamma (a,z)$,  $\Gamma (a,b,z)$,  $\Gamma (a,z,b)$, $Q(a,z)$, 
	$Q(a,b,z)$,  $Q(a,z,b)$, $C_a^{(b)}(z)$, $H_a^{(1)}(z)$, $H_a^{(2)}(z)$,
	$\text{hav}(z)$, $H_a(z)$, $\,_0F_1(;a;z)$, $\,_0\tilde{F}_1(;a;z)$,
	$\,_1F_1(a;b;z)$, $\,_1\tilde{F}_1(a;b;z)$, $\,_2F_1(a,b;c;z)$, $\,_2\tilde{F}_1(a,b;c;z)$,
	$\,_p F_q(a_1,...,a_p;b_1,...,b_q;z)$, $\,_p \tilde{F}_q(a_1,...,a_p;b_1,...,b_q;z)$, $U(a,b,z)$,
	${\rm hav}^{-1}(z)$,
	$\text{bei}_0(z)$,  $\text{bei}_a(z)$, $\text{ber}_0(z)$, $\text{ber}_a(z)$,   $\text{kei}_0(z)$, $\text{kei}_a(z)$,
	$\ker_0(z)$,  $\ker _a(z)$, $L_a(z)$, $L_a^b(z)$,
	$P_{\nu }(z)$,  $ P_a^b(z)$,  $Q_{\nu }(z)$,
	$ Q_a^b(z)$, $\log (z)$, $L_z$, $L_a(z)$, $D_a(z)$,
	$\text{Li}_2(z)$, $\text{Li}_a(z)$,  $\text{Gi}(z)$, $\text{Gi}'(z)$, $\text{Hi}(z)$,
	$\text{Hi}'(z)$, $\sin (z)$, $\text{sinc}(z)$, $ \sinh (z) $, $\text{Shi}(z)$, $\text{Si}(z)$,
	$j_a(z)$, $y_a(z)$, $h_a^{(1)}(z)$, $h_a^{(2)}(z)$,
	$\pmb{H}_{\nu }(z)$, $\pmb{L}_{\nu }(z)$, $\theta (z)$, $\pmb{E}_{\nu }(z)$, $\pmb{E}_{\nu }^a(z)$,
	$M_{a,b}(z)$, $W_{a,b}(z)$ 
	\}
\end{center}

Different combinations, including these functions also sometimes can be presented through Meijer $G$-function. Long time work in this direction allowed to build basic collection of such functions with about 3000 cases and more than 90000 cases were tested.

\subsubsection{Fractional integro-differentiation of Meijer G--functions}

For all above functions one can find different representations of the fractional integro-derivatives searching Wolfram Functions Site \cite{site3}. In particular, definition formula \eqref{FrW} for fractional integro-derivative of the order $\alpha$ of the Meijer $G$-function can be converted to the following representation (see also \cite{site4}):
$$
\frac{d^\alpha}{dz^\alpha}G_{p,q}^{m,n}\left(z\left|
\begin{array}{c}
	a_1,...,a_n,a_{n+1},...,a_p \\
	b_1,...,b_m,b_{m+1},...,b_q \\
\end{array}
\right.\right)
$$
\begin{equation}\label{DG}
	=G_{p+1,q+1}^{m,n+1}\left(z\left|
	\begin{array}{c}
		-\alpha,a_1-\alpha,...,a_n-\alpha,a_{n+1}-\alpha,...,a_p-\alpha \\
		b_1-\alpha,...,b_m-\alpha,0,b_{m+1}-\alpha,...,b_q-\alpha \\
	\end{array}
	\right.\right).
\end{equation}

We can convert the right hand side of this formula to Fox $H$-function if use the following input in \texttt{Wolfram Mathematica}:
\begin{center}
	\texttt{ResourceFunction}["FractionalOrderD"][
	MeijerG[\{Table[Subscript[$a$,$i$], \{i, 1, n\}], 
	Table[Subscript[$a$,$i$], \{i, n + 1, p\}]\}, 
	\{Table[Subscript[$b$,$i$], \{i, 1, m\}], 
	Table[Subscript[$b$,$i$], \{i, m + 1, q\}]\}, 
	z], \{z, $\alpha$\}]
\end{center}

Below we demonstrate more general formula for generalized Meijer $G$-function with argument a $z^r$ and parameter $r_1$:  
$$
\frac{d^\alpha}{dz^\alpha}G_{p,q}^{m,n}\left(az^r,r_1\left|
\begin{array}{c}
	a_1,...,a_n,a_{n+1},...,a_p \\
	b_1,...,b_m,b_{m+1},...,b_q \\
\end{array}
\right.\right)
$$	
\begin{equation}\label{DG1}	=r_1z^{-\alpha}(az^r)^{\alpha/r}H_{p+1,q+1}^{m,n+1}\left(az^r\left|
	\begin{array}{c}
		\{-\alpha,r\},\left\{a_1-\frac{r_1}{r}\alpha,r_1\right\},...,\left\{a_p-\frac{r_1}{r}\alpha,r_1\right\} \\
		\left\{b_1-\frac{r_1}{r}\alpha,r_1\right\},...,\left\{b_q-\frac{r_1}{r}\alpha,r_1 \right\},\{0,r\} \\
	\end{array}
	\right.\right).
\end{equation}

In general (non logarithmic) cases the Meijer $G$-function can be represented through one or finite combination of the series  \eqref{SG} with  $z_0=0$ and $k=0$.

Assume $p\leq q$, 
no two of the bottom parameters $b_j$, $j=1,...,m$, differ by an integer, and $(a_j-b_k)$ 
is not a positive integer when $j=1,2,...,n$ and $k=1,2,...,m$. 
Then double sums of "left"\, residues  for  $G_{p,q}^{m,n}$ can be written as finite sum of generalized hypergeometric functions (see \cite{site22}):
$$
G_{p,q}^{m,n}\left(z \left| 
\begin{array}{c}
	a_1,...,a_n,a_{n+1},...,a_p \\
	b_1,...,b_m,b_{m+1},...,b_q \\
\end{array}\right.
\right)=
$$
\begin{equation}\label{G1}
	=\sum\limits_{k=1}^m A_{p,q,k}^{m,n}\,z^{b_{k}}\,_{p}F_{q-1}\left({\begin{matrix}1+b_{k}-a_{1},...,1+b_{k}-a_{p}\\
			1+b_{k}-b_{1},...*....1+b_{k}-b_{q}\end{matrix}};(-1)^{p-m-n}z\right),
\end{equation}
where $"*"$ indicates that the entry $(1+b_k-b_k)$ is omitted. Also,
$$
A_{p,q,k}^{m,n}={\frac {\prod\limits_{j=1,j\neq k}^{m}\Gamma (b_{j}-b_{k})\prod\limits_{j=1}^{n}\Gamma (1+b_{k}-a_{j})}{\prod\limits_{j=n+1}^{p}\Gamma (a_{j}-b_{k})\prod\limits_{j=m+1}^{q}\Gamma (1+b_{k}-b_{j})}}.
$$
The Function Site \cite{site3} includes the most complete collection of formulas for asymptotics of Meijer $G$-functions.
In more complicated logarithmic cases, the Meijer $G$-function can be represented through a finite combination of the above series \eqref{G1} with $z_0=0$ and $k>0$. It means that we can evaluate Meijer $G$-function, as infinite sums including terms $z^{b}\log^k(z)$ with already uniquely defined basic elementary functions: power and logarithmic.

So if we  have a hypergeometric type function $f(z)$
we first can  write $f(z)$ as $z^b g(z)$  or $z^b \log^k(z) g(z)$  where $g(x)$ is supposed to be simpler than $f(z)$.
Using the series expansion of the function $g(z)=\sum\limits_{n=0}^{\infty } c_n z^n$ and formula \eqref{G1} we can  rewrite $f(z)$ as a single $G$-function. In more complex cases, the function $f(z)$ can be written as a finite sum of $G$-functions.
Next applying  formula \eqref{DG} we can find fractional integral or derivative of  $f(z)$ in the form of the Meijer $G$-function. Then we can write the $G$-function as a simpler function if possible. 

\subsubsection{Fractional integro-differentiation of $e^z$ and $K_\nu(z)$} 

In Wolfram Mathematica there is a function \texttt{MeijerGForm[g(z),z]} which reduces $g(z)$ to the Meijer G--function  as a function of $z$. In order to use function \texttt{MeijerGForm} we should write
\begin{center}
	\texttt{ResourceFunction["MeijerGForm"]}.
\end{center}

For example,
\texttt{ResourceFunction["MeijerGForm"][$e^z$, z]}
gives
$$
e^z=\text{MeijerG}(\{\{\},\{\}\},\{\{0\},\{\}\},-z,1)=G_{0,1}^{1,0}\left(-z,1 \left| 
\begin{array}{c}
	-  \\
	0 \\
\end{array}\right.
\right)=G_{0,1}^{1,0}\left(-z \left| 
\begin{array}{c}
	-  \\
	0 \\
\end{array}\right.
\right)
$$
and \texttt{ResourceFunction["MeijerGForm"][ BesselK[$\nu$, z], z]} gives
$$
K_\nu(z)=\frac{1}{2} \text{MeijerG}\left(\{\{\},\{\}\},\left\{\left\{\frac{\nu }{2},-\frac{\nu
}{2}\right\},\{\}\right\},\frac{z}{2},\frac{1}{2}\right)
=G_{0,2}^{2,0}\left(\frac{z}{2},\frac{1}{2}\left|
\begin{array}{c}
	- \\
	\frac{\nu }{2},-\frac{\nu}{2} \\
\end{array}
\right.\right),
$$
where formula \eqref{Gr} was used.

Let us compare the approach of finding the fractional operator $\frac{d^\alpha }{dz^\alpha}$  through  series expansion and  through   the $G$-function representation.
Direct calculation through  the $G$-function representation gives
\begin{equation}\label{FrE}
	\frac{d^\alpha e^z}{dz^\alpha}=\begin{cases}
		e^z, & \alpha \in \mathbb{Z}\,\, {\rm and}\,\,\alpha \geq 0; \\
		e^z (1-Q(-\alpha ,z)), & \text{in other cases},
	\end{cases}
\end{equation}
where $Q(a,z)=\frac{\Gamma(a,z)}{\Gamma(a)}$ is the regularized incomplete Gamma function,  $\Gamma(a,z)=\int\limits_z^\infty t^{a-1}e^{-t}dt$ is the
incomplete Gamma function.
The exponential function $e^{z}$   has Taylor  series at $z=0$ of the form
$$
e^{z}=\sum _{n=0}^{\infty }{\frac {z^{n}}{n!}}.
$$
Since $n\in\mathbb{N}\cup\{0\}$ using the third line in \eqref{powerG} we obtain
\begin{align*}
	\frac{d^\alpha e^z}{dz^\alpha}&=\sum _{n=0}^{\infty }{\frac {1}{n!}}\frac{d^\alpha z^{n}}{dz^\alpha}\\
	&=\sum _{n=0}^{\infty }{\frac {1}{n!}}\frac{\Gamma (n+1) z^{n-\alpha }}{\Gamma (n-\alpha +1)}=\sum _{n=0}^{\infty } \frac{z^{n-\alpha }}{\Gamma (n-\alpha +1)}\\
	&=e^z\left(1+\frac{\alpha\Gamma (-\alpha ,z)}{\Gamma (1-\alpha )}\right).
\end{align*}
For $\alpha=0,1,2,...$ we have $\frac{\alpha\Gamma (-\alpha ,z)}{\Gamma (1-\alpha )}=0$ and $$\frac{\alpha\Gamma (-\alpha ,z)}{\Gamma (1-\alpha )}=\frac{\alpha\Gamma (-\alpha ,z)}{-\alpha\Gamma (-\alpha )}=-\frac{\Gamma (-\alpha ,z)}{\Gamma (-\alpha )}=-Q(a,z),$$
therefore, the result coincides with  formula \eqref{FrE}.

For $K_0(z)$ we have
\begin{equation}\label{FrK01}
	\frac{d^\alpha }{dz^\alpha}K_0(z)=\frac{1}{2} G_{2,4}^{2,2}\left(\frac{z}{2},\frac{1}{2}\left|
	\begin{array}{c}
		\frac{1-\alpha }{2},-\frac{\alpha }{2} \\
		-\frac{\alpha }{2},-\frac{\alpha }{2},0,\frac{1}{2} \\
	\end{array}
	\right.\right).
\end{equation}
On the other hand by \eqref{H003}
$$
\frac{d^\alpha }{dz^\alpha} K_0(z)=\frac{d^\alpha }{dz^\alpha}\left(  -\left( \log({z})-\log(2) +\gamma\right)  \sum\limits_{n=0}^\infty\frac{1}{(n!)^2}\left(\frac{z}{2}\right)^{2n} +\sum\limits_{n=1}^\infty  \frac{H_n}{(n!)^2}\left({\frac {z}{2}}\right)^{2n}\right) ,
$$
where $\gamma$  is Euler--Mascheroni constant \eqref{EM}. Applying \eqref{powerG} and \eqref{DLog} we obtain
\begin{align}\label{FrK02}
	\frac{d^\alpha }{dz^\alpha} K_0(z)&\nonumber\\
	&	=\frac{1}{2}\sum\limits_{k=0}^\infty\frac{2^{1-2 k} (2 k)! \psi(k+1) }{(k!)^2 \Gamma (2 k-\alpha +1)}z^{2 k-\alpha }\nonumber\\
	&-\sum\limits_{k=0}^\infty\frac{   (2 k)!  \left(H_{2 k}-H_{2 k-\alpha }+\log \left(\frac{z}{2}\right)\right)}{2^{2 k}(k!)^2
		\Gamma (2 k-\alpha +1)}z^{2 k-\alpha }.
\end{align}
Using Mathematica one  can check that \eqref{FrK01} equals to \eqref{FrK02}.

\subsubsection{Riemann-Liouville integral of generalized Meijer G--function}

Let us now  evaluate the following "main"\, Riemann-Liouville integral
\begin{equation}\label{GrInt}
	\frac{1}{\Gamma(\beta)}
	\int\limits_0^z (z-\tau)^{\beta-1}\tau^{\alpha-1}
	G_{p,q}^{m,n}\left(w\tau^g,r \left| 
	\begin{array}{c}
		a_1,...,a_n,a_{n+1},...,a_p \\
		b_1,...,b_m,b_{m+1},...,b_q \\
	\end{array}\right.
	\right)d\tau
\end{equation}
which defines fractional order $\beta$ 
integral of generalized Meijer $G$-function \eqref{Gr}, 
multiplied by $\tau^{\alpha-1}$, with  $r\in \mathbb{R}$, $r\neq 0$, $m\in \mathbb{Z}$, $m\geq 0$, $n\in \mathbb{Z}$, $n\geq 0$, $p\in
\mathbb{Z}$, $p\geq 0$, $q\in \mathbb{Z}$, $q\geq 0$, $m\leq q$, $n\leq p$.
In formula  \eqref{Gr} for simplicity we use vertical contour 
$\mathscr{L}=\{\gamma-i\infty,\gamma+i\infty\}$.

The integral \eqref{GrInt} with particular parameters $g=\frac{\ell}{k}\in\mathbb{Q}$ and 
$r=1$ has value, shown at \cite{MarichevMellin}, p. 535, formula 3.36.2.1. 
Below we derive representations of this integral through Fox $H$-function or Meijer $G$-function and  demonstrate in
detail how to build the corresponding sets of conditions for convergence of this integral.

After substitution of above definition of Meijer $G$-function to Riemann-Liouville integral and  
changing the order of integration with respect to
$s$ and $\tau$
we come to the following chain operations where we met and evaluated internal  integral
\begin{align}\label{LF}
	I&=\frac{1}{\Gamma(\beta)}
	\int\limits_0^z (z-\tau)^{\beta-1}\tau^{\alpha-1}
	G_{p,q}^{m,n}\left(w\tau^g,r \left| 
	\begin{array}{c}
		a_1,...,a_n,a_{n+1},...,a_p \\
		b_1,...,b_m,b_{m+1},...,b_q \\
	\end{array}\right.
	\right)d\tau\nonumber\\
	&=\frac{1}{2 \pi i \Gamma (\beta )}\int\limits_0^z (z-\tau)^{\beta-1}\tau^{\alpha-1}\times\nonumber\\
	&\times\left(\int\limits_{\gamma-i\infty}^{\gamma+i\infty}\frac{\prod\limits_{k=1}^{m}\Gamma (b_{k}+s)\prod\limits_{k=1}^{n}\Gamma (1-a_{k}-s)}{\prod\limits_{k=n+1}^{p}\Gamma (a_{k}+s)\prod\limits_{k=m+1}^{q}\Gamma (1-b_{k}-s)}\,(w\tau^g)^{-\frac{s}{r}}\,ds\right) d\tau\nonumber\\
	&=\frac{1}{2 \pi i \Gamma (\beta )}\int\limits_{\gamma-i\infty}^{\gamma+i\infty}
	w^{-\frac{s}{r}} \frac{\prod\limits_{k=1}^{m}\Gamma (b_{k}+s)\prod\limits_{k=1}^{n}\Gamma (1-a_{k}-s)}{\prod\limits_{k=n+1}^{p}\Gamma (a_{k}+s)\prod\limits_{k=m+1}^{q}\Gamma (1-b_{k}-s)}\times\nonumber\\
	&\times \left( \int\limits_0^z (z-\tau)^{\beta-1}\tau^{\alpha-g\frac{s}{r}-1} d\tau\right)\,ds\nonumber\\
	&	=\frac{z^{\alpha+\beta-1}}{2 \pi i }\int\limits_{\gamma-i\infty}^{\gamma+i\infty}
	\frac{ \prod\limits_{k=1}^m \Gamma \left(b_k+s\right)\prod\limits_{k=1}^n \Gamma
		\left(1-a_k-s\right)\Gamma \left(\alpha -\frac{g }{r}s\right)}{\prod\limits_{k=n+1}^p \Gamma \left(a_k+s\right)\prod\limits_{k=m+1}^q \Gamma \left(1-b_k-s\right)\Gamma \left(\alpha +\beta-\frac{g}{r}s
		\right)} \times\nonumber\\
	&\times \left(w^{\frac{1}{r}} z^{\frac{g}{r}}\right)^{-s}ds.
\end{align}
The last Mellin-Barnes integral \eqref{LF} can be written through the corresponding Fox 
$H$-function by definition of Fox $H$-function. 
As a result we have the following representation of the Riemann-Liouville fractional order $\beta$ 
integral of the generalized Meijer $G$-function through Fox $H$-function
$$
\frac{1}{\Gamma(\beta)}
\int\limits_0^z (z-\tau)^{\beta-1}\tau^{\alpha-1}
G_{p,q}^{m,n}\left(w\tau^g,r \left| 
\begin{array}{c}
	a_1,...,a_n,a_{n+1},...,a_p \\
	b_1,...,b_m,b_{m+1},...,b_q \\
\end{array}\right.
\right)d\tau
$$
\begin{equation}\label{LF1}
	=z^{\alpha +\beta -1}  H_{p,q}^{\,m,n}\left[w^{\frac{1}{r}} z^{\frac{g}{r}}\left|{\begin{matrix}\left(1-\alpha,\frac{g}{r}\right),(a_{1},1),...,(a_{n},1),(a_{n+1},1),...,(a_p,1)\\
			(b_{1},1),...,(b_{m},1),(b_{m+1},1),...,(b_q,1),\left(1-\alpha-\beta,\frac{g}{r}\right) \end{matrix}}\right.\right].
\end{equation}
We should mention that operations like $(w\tau^g)^{-\frac{s}{r}}\rightarrow \left(w^{\frac{1}{r}}z^{\frac{g}{r}} \right)^{-s} $ with complex variables demand special conditions for their correctness.

The   formula \eqref{LF} with products of Gamma functions includes gammas with variable of integration $s$ which has coefficients $+1$ or $-1$ or $-\frac{g}{r}$. 
Definition of Meijer $G$-function does not include coefficients not equal to $+1$ or $-1$ 
but after some special transformations above Mellin-Barnes integral can be re-written as Meijer $G$-function.

Let consider the Mellin-Barnes integral from  formula  \eqref{LF}, which we denote as $MB$:
\begin{align*}
	MB&=\frac{z^{\alpha+\beta-1}}{2 \pi i}\int\limits_{\gamma-i\infty}^{\gamma+i\infty}
	\frac{\Gamma \left(\alpha -\frac{g }{r}s\right)}{\Gamma \left(\alpha +\beta-\frac{g}{r}s
		\right)}\cdot\\
	&\cdot\frac{ \prod\limits_{k=1}^m \Gamma \left(b_k+s\right)}{\prod\limits_{k=n+1}^p \Gamma \left(a_k+s\right)} \cdot \frac{\prod\limits_{k=1}^n \Gamma
		\left(1-a_k-s\right)}{\prod\limits_{k=m+1}^q \Gamma \left(1-b_k-s\right)} \left(w^{\frac{1}{r}} z^{\frac{g}{r}}\right)^{-s}ds.
\end{align*}
Suppose that $\frac{g}{r}$ is the rational value with positive integers $g$ and $r$ (that do not have common factors).
Making change of variable $s= r\zeta$ we obtain
\begin{align*}
	MB& =\frac{z^{\alpha+\beta-1}}{2 \pi i}\int\limits_{\delta-i\infty}^{\delta+i\infty}
	\frac{\Gamma \left(\alpha -g\zeta\right)}{\Gamma \left(\alpha +\beta-g\zeta
		\right)}\cdot\\
	&\cdot\frac{ \prod\limits_{k=1}^m \Gamma \left(b_k+r\zeta\right)}{\prod\limits_{k=n+1}^p \Gamma \left(a_k+r\zeta\right)} \cdot \frac{\prod\limits_{k=1}^n \Gamma
		\left(1-a_k-r\zeta\right)}{\prod\limits_{k=m+1}^q \Gamma \left(1-b_k-r\zeta\right)} \left(w^{\frac{1}{r}} z^{\frac{g}{r}}\right)^{-r\zeta} rd\zeta,\quad \delta=\frac{\gamma}{r}.
\end{align*}

Since $r,g\in\mathbb{N}$ for each Gamma function in the previous Mellin-Barnes integral,
we can apply known  Gauss multiplication formulas  of the forms (see \cite{site5})
\begin{equation}\label{GF01}
	\Gamma (b+\zeta  r)=(2 \pi )^{\frac{1-r}{2}} r^{b+\zeta  r-\frac{1}{2}} \prod\limits_{j=0}^{r-1} \Gamma
	\left(\frac{b+j}{r}+\zeta \right),\qquad r\in\mathbb{N},
\end{equation}
\begin{equation}\label{GF02}
	\Gamma (\alpha -\zeta  g)=(2 \pi )^{\frac{1-g}{2}} g^{\alpha -\zeta  g-\frac{1}{2}}
	\prod\limits_{j=0}^{g-1} \Gamma
	\left(\frac{\alpha+j}{g}-\zeta \right),\qquad g\in\mathbb{N}
\end{equation}
for each Gamma function.
We get
\begin{align*}
	\mathscr{A}(\alpha,\beta,r)&=\frac{\Gamma (\alpha -g \zeta )}{\Gamma (\alpha +\beta -g \zeta )}\cdot
	\frac{\prod\limits_{k=1}^m \Gamma \left(b_k+r \zeta\right)}{\prod\limits_{k=n+1}^p \Gamma \left(a_k+r \zeta\right)}\cdot\\
	&\cdot \frac{\prod\limits_{k=1}^n \Gamma \left(1-a_k-r \zeta\right)}{\prod\limits_{k=m+1}^q \Gamma \left(1 -b_k-r \zeta\right)}\cdot 
	\left(w^{\frac{1}{r}} z^{\frac{g}{r}}\right)^{-r\zeta}r =\\
	&=\frac{(2 \pi )^{\frac{1-g}{2}}}{(2\pi)^\frac{1-g}{2}}\cdot \frac{g^{\alpha-g\zeta-\frac{1}{2}}}{g^{\alpha+\beta-g\zeta-\frac{1}{2}}}\cdot \frac{\prod\limits _{j=0}^{g-1} \Gamma \left(\frac{\alpha+j}{g}-\zeta \right)}{\prod\limits_{j=0}^{g-1} \Gamma
		\left(\frac{\alpha +\beta+j}{g}-\zeta \right)}\cdot \\
	&	
	\cdot\frac{\left(\prod\limits_{k=1}^n (2 \pi )^{\frac{1-r}{2}} r^{-a_k-r\zeta+\frac{1}{2}}\right) \prod\limits
		_{j=0}^{r-1} \Gamma \left(\frac{j-a_k+1}{r}-\zeta \right)}{ \left(\prod\limits_{k=m+1}^q (2 \pi )^{\frac{1-r}{2}}
		r^{-b_k-r\zeta+\frac{1}{2}}\right)\prod\limits_{j=0}^{r-1} \Gamma
		\left(\frac{j-b_k+1}{r}-\zeta \right)}\cdot\\
	&
	\cdot \frac{ \left(\prod\limits_{k=1}^m (2 \pi )^{\frac{1-r}{2}} r^{b_k+r\zeta 
			-\frac{1}{2}}\right) \prod\limits_{j=0}^{r-1} \Gamma \left(\frac{b_k+j}{r}+\zeta\right)}{ \left(\prod\limits_{k=n+1}^p (2 \pi
		)^{\frac{1-r}{2}} r^{a_k+r\zeta-\frac{1}{2}}\right)\prod\limits_{j=0}^{r-1} \Gamma \left(\frac{a_k+j}{r}+\zeta\right)} \left(w z^g\right)^{-\zeta }r.
\end{align*}
Simplifying we obtain
\begin{align*}
	\mathscr{A}(\alpha,\beta,r)&=\\
	&=
	\frac{g^{-\beta } \prod\limits_{j=0}^{g-1} \Gamma \left(\frac{\alpha+j}{g}-\zeta \right)}{\prod\limits
		_{j=0}^{g-1} \Gamma \left(\frac{\alpha +\beta+j}{g}-\zeta \right)}\cdot\\
	&\cdot\frac{\prod\limits _{k=1}^n (2 \pi )^{\frac{1-r}{2}} r^{-a_k-r\zeta+\frac{1}{2}} \prod\limits_{j=0}^{r-1}
		\Gamma \left(\frac{1-a_k+j}{r}-\zeta \right)}{\prod\limits _{k=m+1}^q (2 \pi )^{\frac{1-r}{2}}
		r^{-b_k-r\zeta+\frac{1}{2}} \prod\limits_{j=0}^{r-1} \Gamma \left(\frac{1-b_k+j}{r}-\zeta
		\right)}\cdot\\
	&
	\cdot \frac{\prod\limits_{k=1}^m (2 \pi )^{\frac{1-r}{2}} r^{b_k+r\zeta -\frac{1}{2}} \prod\limits_{j=0}^{r-1}
		\Gamma \left(\frac{b_k+j}{r}+\zeta \right)}{\prod\limits_{k=n+1}^p (2 \pi )^{\frac{1-r}{2}}
		r^{a_k+r\zeta -\frac{1}{2}} \prod\limits_{j=0}^{r-1} \Gamma \left(\frac{a_k+j}{r}+\zeta\right)} \left(w z^g\right)^{-\zeta }r=\\
	&
	=\frac{\prod\limits_{k=1}^n (2 \pi )^{\frac{1-r}{2}} r^{-a_k-r\zeta+\frac{1}{2}}}
	{\prod\limits_{k=m+1}^q (2
		\pi )^{\frac{1-r}{2}} r^{-b_k-r\zeta+\frac{1}{2}}}\cdot \frac{\prod\limits_{k=1}^m (2 \pi )^{\frac{1-r}{2}} r^{b_k+r\zeta-\frac{1}{2}}}{\prod\limits_{k=n+1}^p (2 \pi
		)^{\frac{1-r}{2}} r^{a_k+r\zeta -\frac{1}{2}}}\cdot
	\frac{\prod\limits_{j=0}^{g-1} \Gamma \left(\frac{\alpha+j}{g}-\zeta \right)}
	{\prod\limits_{j=0}^{g-1} \Gamma
		\left(\frac{\alpha +\beta+j}{g}-\zeta \right)}\cdot\\
	&\cdot \frac{\prod\limits_{k=1}^n \prod\limits_{j=0}^{r-1} \Gamma \left(\frac{1-a_k+j}{r}-\zeta \right)}{\prod\limits_{k=m+1}^q \prod\limits_{j=0}^{r-1} \Gamma \left(\frac{1-b_k+j}{r}-\zeta \right)} \cdot
	\frac{ \prod\limits_{k=1}^m \prod\limits_{j=0}^{r-1} \Gamma
		\left(\frac{b_k+j}{r}+\zeta\right)}{\prod\limits_{k=n+1}^p \prod\limits_{j=0}^{r-1} \Gamma \left(\frac{a_k+j}{r}+\zeta\right)} g^{-\beta } \left(w z^g\right)^{-\zeta } r.
\end{align*}
Collecting powers we get
\begin{align*}
	\mathscr{A}(\alpha,\beta,r)&= \frac{\left((2 \pi )^{\frac{1-r}{2}} r^{\frac{1}{2}-r\zeta}\right)^{n}}{\left((2 \pi )^{\frac{1-r}{2}} r^{\frac{1}{2}-r\zeta}\right)^{q-m}}\cdot \frac{\left((2 \pi )^{\frac{1-r}{2}} r^{\zeta  r-\frac{1}{2}}\right)^{m}}{\left((2 \pi )^{\frac{1-r}{2}} r^{\zeta  r-\frac{1}{2}}\right)^{p-n}}\cdot\\
	&\cdot \frac{\prod\limits_{k=1}^n r^{-a_k}}{\prod\limits_{k=m+1}^q r^{-b_k}}\cdot \frac{\prod\limits_{k=1}^m r^{b_k}}{\prod\limits_{k=n+1}^p r^{a_k}}\cdot\frac{\prod\limits_{j=0}^{g-1} \Gamma \left(\frac{\alpha+j }{g}-\zeta \right)}{\prod\limits_{j=0}^{g-1} \Gamma
		\left(\frac{\alpha +\beta+j }{g}-\zeta \right)}\cdot\\
	&
	\cdot \frac{\prod\limits_{k=1}^n \prod\limits _{j=0}^{r-1} \Gamma \left(\frac{1-a_k+j}{r}-\zeta \right)}{\prod\limits
		_{k=m+1}^q \prod\limits _{j=0}^{r-1} \Gamma \left(\frac{1-b_k+j}{r}-\zeta \right)}\cdot
	\frac{\prod\limits_{k=1}^m \prod\limits _{j=0}^{r-1} \Gamma \left(\frac{b_k+j}{r}+\zeta\right)}{\prod\limits
		_{k=n+1}^p \prod\limits_{j=0}^{r-1} \Gamma \left(\frac{a_k+j}{r}+\zeta\right)}g^{-\beta } \left(w z^g\right)^{-\zeta } r=\\
	&
	=\frac{\prod\limits_{j=0}^{g-1} \Gamma \left(\frac{\alpha+j}{g}-\zeta \right)}
	{\prod\limits_{j=0}^{g-1} \Gamma
		\left(\frac{\alpha +\beta+j}{g}-\zeta \right)}\cdot\frac{\prod\limits_{k=1}^n 
		\prod\limits_{j=0}^{r-1} \Gamma \left(\frac{1-a_k+j}{r}-\zeta \right)}{\prod\limits
		_{k=m+1}^q \prod\limits_{j=0}^{r-1} \Gamma \left(\frac{1-b_k+j}{r}-\zeta \right)}\cdot\\
	&
	\cdot\frac{\prod\limits_{k=1}^m \prod\limits_{j=0}^{r-1} \Gamma \left(\frac{b_k+j}{r}+\zeta\right)}{\prod\limits
		_{k=n+1}^p \prod\limits_{j=0}^{r-1} \Gamma \left(\frac{a_k+j}{r}+\zeta\right)}\cdot \frac{ r^{\sum\limits_{k=1}^q
			b_k-\sum\limits _{k=1}^p a_k+\frac{p-q}{2}+1}}{(2 \pi )^{(r-1)\left(m+n-\frac{p+q}{2}\right)}g^{\beta }} \cdot \left(\frac{w z^g}{ r^{r (q-p)}}\right)^{-\zeta }.
\end{align*}
In above transformations we changed $\left(w^{\frac{1}{r}} z^{\frac{g}{r}}\right)^{-r\zeta}$ to $(w z^g)^{-\zeta}$, 
which is allowed under some restrictions. 
But in the result we arrived at products of Gamma functions with coefficients $+1$ or $-1$ at the variable of integration $\zeta$.

This result allows us to write the following representation of Riemann-Liouville fractional order $\beta$
integral of the generalized Meijer $G$-function through Meijer $G$-function in the case when $\frac{g}{r}$ 
is the rational value with positive integers $g$ and $r$:
$$
\frac{1}{\Gamma(\beta)}
\int\limits_0^z (z-\tau)^{\beta-1}\tau^{\alpha-1}
G_{p,q}^{m,n}\left(w\tau^g,r \left| 
\begin{array}{c}
	a_1,...,a_n,a_{n+1},...,a_p \\
	b_1,...,b_m,b_{m+1},...,b_q \\
\end{array}\right.
\right)d\tau
$$
\begin{equation}\label{FrIntG}
	=z^{\alpha+\beta-1}\frac{r^{\sum\limits_{k=1}^q b_k-\sum\limits_{i=1}^p
			a_i+\frac{p-q}{2}+1}}{ (2 \pi )^{(r-1) \left(m+n-\frac{p+q}{2}\right)} g^{\beta }} G_{rp+g,rq+g}^{rm,rn+g}\left(\frac{wz^g}{r^{r(q-p)}} \left| 
	\begin{array}{c}
		\Delta_1,\Delta_2 \\
		\Delta_3, \Delta_4\\
	\end{array}\right.\right),
\end{equation}
where 
$$\Delta_1= \frac{1-\alpha }{g},...,\frac{g-\alpha}{g},\frac{a_1}{r},...,\frac{a_1+r-1}{r},...,\frac{a_n}{r},...,\frac{a_n+r-1}{r},$$
$$
\Delta_2= \frac{a_{n+1}}{r},...,\frac{a_{n+1}+r-1}{r},...,\frac{a_p}{r},...,\frac{a_p+r-1}{r}, 
$$
$$
\Delta_3= \frac{b_1}{r},...,\frac{b_1+r-1}{r},...,\frac{b_m}{r},...,\frac{b_m+r-1}{r},
$$
$$
\Delta_4=\frac{b_{m+1}}{r},...,\frac{b_{m+1}+r-1}{r},...,\frac{b_q}{r},...,\frac{b_q+r-1}{r},\frac{1-\alpha-\beta}{g},...,\frac{g-\alpha-\beta}{g}.
$$

\subsubsection{Big O representations for asymptotics of Meijer $G$-functions}

If we use "Big O notation"\, \cite{BigO} for description of asymptotic behaviour of Meijer $G$-function, we get the following generic formulas, describing Big O representation of Meijer $G$-function near its two or three singular points (see Wolfram Functions Site \cite{site} and \cite{MarichevMellin}, p. 571):
\begin{align}\label{G33}
	G_{p,q}^{m,n}\left(z \left| 
	\begin{array}{c}
		a_1,...,a_n,a_{n+1},...,a_p \\
		b_1,...,b_m,b_{m+1},...,b_q \\
	\end{array}\right.\right)\nonumber\\
	\leftrightarrow \left\{ \begin{array}{ll}
		\sum\limits_{k=1}^m z^{b_k} & \mbox{if $p\leq q$};\\
		\sum\limits_{k=1}^m z^{b_k}+z^{\chi } e^{\frac{(-1)^{q-m-n}}{z}} & \mbox{if $p=q+1$};\\
		\sum\limits_{k=1}^m z^{b_k}+z^{\chi } \cos \left(2 \sqrt{\frac{(-1)^{q-m-n-1}}{z}}\right) & \mbox{if $p=q+2$};\\
		\sum\limits_{k=1}^m z^{b_k}+z^{\chi } e^{(p-q) (-z)^{\frac{1}{q-p}}} & \mbox{if $p\geq q+3$},\end{array} \right.
\end{align}
in \eqref{G33} $|z|\rightarrow 0$,
\begin{equation}\label{chi}
	\chi =\frac{1}{{q-p}}\left( \sum\limits _{j=1}^q b_j-\sum\limits_{j=1}^p a_j+\frac{p-q+1}{2} \right).
\end{equation}
When $z\to (-1)^{m+n-p}$ we get
\begin{align}\label{GAssim}
	G_{p,q}^{m,n}\left(z \left| 
	\begin{array}{c}
		a_1,...,a_n,a_{n+1},...,a_p \\
		b_1,...,b_m,b_{m+1},...,b_q \\
	\end{array}\right.\right)\nonumber \\ \leftrightarrow
	\left\{
	\begin{array}{ll}
		1+\left(1- (-1)^{p-m-n}z\right)^{\psi_p} & \mbox{if $q=p$ and $\psi _p\neq 0$}; \\
		1+\log \left(1- (-1)^{p-m-n}z\right)& \mbox{if $q=p$ and $\psi _p=0$};  \\
		1 & \mbox{if $q\neq p$}, 
	\end{array} \right.
\end{align}
where 
\begin{equation}\label{psi}
	\psi _p=\sum _{j=1}^p \left(a_j-b_j\right)-1.
\end{equation}
Finally, for $|z|\rightarrow \infty$
\begin{align}\label{GAssim2}
	G_{p,q}^{m,n}\left(z \left| 
	\begin{array}{c}
		a_1,...,a_n,a_{n+1},...,a_p \\
		b_1,...,b_m,b_{m+1},...,b_q \\
	\end{array}\right.\right)\nonumber \\
	\leftrightarrow \left\{ \begin{array}{ll}
		\sum\limits_{k=1}^n z^{a_k-1} & \mbox{if $q\leq p$};\\
		\sum\limits_{k=1}^n z^{a_k-1}+z^{\chi } e^{{(-1)^{p-m-n}}{z}} & \mbox{if $q=p+1$};\\
		\sum\limits_{k=1}^n z^{a_k-1}+z^{\chi } \cos \left(2 \sqrt{{(-1)^{p-m-n-1}}{z}}\right) & \mbox{if $q=p+2$};\\
		\sum\limits_{k=1}^n z^{a_k-1}+z^{\chi } e^{(q-p) (-z)^{\frac{1}{q-p}}} & \mbox{if $q\geq p+3$},\end{array} \right.
\end{align}
where $\chi$ is defined by \eqref{chi}.

\subsubsection{Conditions of convergence for Riemann-Liouville integrals with Meijer $G$-functions}

Above relations include power functions like $z^{b_k}$ or $z^\chi$ or $z^{a_k-1}$ or $1$ or $\log(1- (-1)^{p-m-n}z)$ 
which match "Big O terms"\, \cite{BigO} and used here instead of $O(z^{b_k})$ or $O(z^\chi)$ or $O(z^{a_k-1})$ or $O(1)$ or $O(\log(1- (-1)^{p-m-n}z))$ inside of asymptotic expansions. It allows for a very simple
way to establish conditions of convergence integrals involving Meijer $G$-functions. For example, the classical Riemann-Liouville integral
\begin{equation}\label{GINT}
	\frac{1}{\Gamma(\beta)}
	\int\limits_0^z (z-\tau)^{\beta-1}\tau^{\alpha-1}
	G_{p,q}^{m,n}\left(w\tau^g,r \left| 
	\begin{array}{c}
		a_1,...,a_n,a_{n+1},...,a_p \\
		b_1,...,b_m,b_{m+1},...,b_q \\
	\end{array}\right.
	\right)d\tau
\end{equation}
has Meijer $G$-function $G_{p,q}^{m,n}[w \tau^g,r|...]$ with parameter  $r$, which under some conditions can be written as classical  $G$-function without parameter $r=1$ as $G_{p,q}^{m,n}[w^{\frac{1}{r}} \tau^{\frac{g}{r}}|...]$. 
The interval of integration here is finite andwe do not need
conditions for convergence at infinity if $\frac{g}{r}>0$ (then $\tau^\frac{g}{r}\rightarrow \infty$ for $\tau\rightarrow\infty$). But for $q<p$ and $\frac{g}{r}>0$ the point 
$\tau=0$ becomes essentially singular and this integral will converge at zero if and only if the following "Big O-equivalent"\, integrals will converge:
\begin{equation}\label{first}
	\int\limits_0^z \tau^{\alpha-1}\sum\limits_{k=1}^m z^{b_k}d\tau,\qquad p\leq q,
\end{equation}
\begin{equation}\label{second}
	\int\limits_0^z \left( \tau^{\alpha-1}\sum\limits_{k=1}^m z^{b_k}+\tau^{\alpha-1}z^\chi e^{\frac{(-1)^{q-m-n}}{z}}\right) d\tau,\qquad p=q+1,
\end{equation}
\begin{equation}\label{third}
	\int\limits_0^z \left( \tau^{\alpha-1}\sum\limits_{k=1}^m z^{b_k}+\tau^{\alpha-1}z^\chi \cos\left(2\sqrt{{\frac{(-1)^{q-m-n}}{z}}}\right) \right) d\tau,\qquad p=q+2,
\end{equation}
\begin{equation}\label{last}
	\int\limits_0^z \left( \tau^{\alpha-1}\sum\limits_{k=1}^m z^{b_k}+\tau^{\alpha-1}z^\chi e^{(p-q)(-z)^{\frac{1}{q-p}}} \right) d\tau,\qquad p\geq q+3,
\end{equation}
where  $z=w^{\frac{1}{r}}\tau^{\frac{g}{r}}$ and  $\frac{g}{r}>0$, $\chi$ is defined by \eqref{chi}. 
The convergence of the integral \eqref{first} takes place when the  integral \eqref{last} in the
above chain converges.
It can happened under condition
\begin{equation}\label{min}
	{\rm M}={\rm Re}(\alpha)+\frac{g}{r}\min\limits_{1\leq k\leq m}\{{\rm Re}(b_k)\}>0.
\end{equation}
Evidently, that this condition appears in other three situation with $p>q$. 
Let 
$$\theta (\delta ) = \left\{ 
\begin{array}{c}
	0,\,\delta<0; \\
	1,\,\delta\geq0. \\
\end{array}\right.
$$
There we see three additional integrals with exponential and cosine functions, for which one can write and evaluate the following "convergent equivalents"\,:
$$
\int\limits_{0}^1 x^\gamma e^{ax^\delta}dx=\frac{(-a)^{-\frac{\gamma +1}{\delta }} }{\delta }\left(\theta (\delta ) \Gamma \left(\frac{\gamma +1}{\delta }\right)-\Gamma
\left(\frac{\gamma +1}{\delta },-a\right)\right),
$$
where ${\rm Re}(\gamma)>-1$ when $\delta>0$ and ${\rm Re}(a)\leq 0$ when $\delta<0$;
$$
\int\limits_{0}^1 x^\gamma \cos({ax^\delta})dx=\frac{1}{\gamma +1}\, _1F_2\left(\frac{\gamma+1}{2 \delta };\frac{1}{2},\frac{\gamma+1 }{2 \delta }+1;-\frac{a^2}{4}\right)-
$$
$$
-\frac{\theta (-\delta )}{\delta } \cos \left(\frac{\pi  (\gamma +1)}{2 \delta }\right) \Gamma
\left(\frac{\gamma +1}{\delta }\right) | a| ^{-\frac{\gamma +1}{\delta }},
$$
where $a\in \mathbb{R}$, $\delta \in \mathbb{R}$ and ${\rm Re}\left(\frac{\gamma +1}{\delta }\right)<1$.

If we apply above conditions of convergence to mentioned three integrals, re-written in the form
$$
\int\limits_0^1 \tau^{\alpha +\frac{g}{r}\chi-1} e^{(-1)^{q-m-n}w^{-\frac{1}{r}} \tau ^{-\frac{g}{r}} }d\tau,
$$
$$
\int\limits_0^1 \tau^{\alpha +\frac{g}{r}\chi-1} \cos \left(2 \tau ^{\frac{g}{2 r}} \sqrt{(-1)^{q-m-n-1}w^{-\frac{1}{r}} }\right)d\tau,
$$
$$
\int\limits_0^1 \tau ^{\alpha +\frac{g}{r}\chi-1} e^{(p-q) \left(-w^{\frac{1}{r}}\right)^{\frac{1}{q-p}} \tau ^{\frac{g}{r
			(q-p)}}} d\tau
$$
we come to the following three groups of conditions
\begin{equation}\label{cond01}
	\left\{
	\begin{array}{cc}
		{\rm Re}\left(\alpha +\frac{g  }{r}\chi\right)>0, & g r<0; \\
		{\rm Re}\left(w^{-\frac{1}{r}} (-1)^{q-m-n}\right)\leq 0, & g r>0, \\
	\end{array}
	\right.
\end{equation}
\begin{equation}\label{cond02}
	\sqrt{w^{-\frac{1}{r}} (-1)^{q-m-n-1}}\in \mathbb{R}\quad\text{and}\quad\frac{g}{r}\in \mathbb{R}\quad\text{and}\quad 2 {\rm Re}\left(\frac{\alpha 
		r}{g}+\chi \right)<1,
\end{equation}
\begin{equation}\label{cond03}
	\left\{
	\begin{array}{cc}
		{\rm Re}\left(\alpha +\frac{g }{r}\chi\right)>0, & \frac{g}{r (q-p)}>0; \\
		{\rm Re}\left((p-q) \left(-w^{\frac{1}{r}}\right)^{\frac{1}{q-p}}\right)\leq 0, & \frac{g}{r (q-p)}<0. \\
	\end{array}
	\right.
\end{equation}
After incorporating these conditions into "Big O-equivalent"\, integrals above we come to the following conditions of convergence of the initial Riemann-Liouville integral at the point $\tau=0$:
\begin{equation}\label{cond04}
	\left\{ \begin{array}{ll}
		{\rm M}>0,\quad &  \mbox{$p\leq q$};\\
		{\rm M}>0,\quad \left\{
		\begin{array}{cc}
			{\rm Re}\left(\alpha +\frac{g}{r}\chi \right)>0, & \frac{g}{r}<0; \\
			{\rm Re}\left((-1)^{q-m-n}w^{-\frac{1}{r}} \right)\leq 0, & \frac{g}{r}>0, \\
		\end{array}
		\right.\qquad &  \mbox{$p=q+1$};\\
		{\rm M}>0,\quad
		\sqrt{w^{-\frac{1}{r}} (-1)^{q-m-n-1}}\in \mathbb{R},\,\, \frac{g}{r}\in \mathbb{R},\,\, {\rm Re}\left(\frac{\alpha  r}{g}+\chi
		\right)<\frac{1}{2},\,&  \mbox{$p=q+2$};\\
		{\rm M}>0,\quad
		\left\{
		\begin{array}{cc}
			{\rm Re}\left(\alpha +\frac{g }{r}\chi \right)>0, & \frac{g}{r (q-p)}>0; \\
			{\rm Re}\left((p-q) \left(-w^{\frac{1}{r}}\right)^{\frac{1}{q-p}}\right)\leq 0, & \frac{g}{r (q-p)}<0, \\
		\end{array}
		\right.\quad &  \mbox{$p\geq q+3$},
	\end{array} \right. 
\end{equation}
with $\frac{g}{r}>0$, $\chi$ is given by \eqref{chi}.

Evidently, that we need to add the restriction ${\rm Re}(\beta)>0$
for convergence of the integral at point $\tau=z$ 
and maybe we can have third singular point of Meijer $G$-function  $w^{\frac{1}{r}}\tau^{\frac{g}{r}}=(-1)^{m+n-p}$, 
arising for $q=p$.  

Let
\begin{equation}\label{tau}
	\tau_0=((-1)^{m+n-p} w^{\frac{1}{r}})^{\frac{r}{g}}. 	
\end{equation}
Then under conditions $\tau_0\in\mathbb{R}$
and $0<\tau_0<z$ 
we should add the restriction 
\begin{equation}\label{cond05}
	{\rm Re}\left(\sum\limits_{j=1}^p(a_j-b_j)\right)>0
\end{equation}
(it is ${\rm Re}(\psi_p+1)>0$ for $q=p$) 
which becomes weaker as
\begin{equation}\label{cond06}
	{\rm Re}\left(\sum\limits_{j=1}^p(a_j-b_j)+\beta-1\right)>0
\end{equation}
if point $\tau_0$ coincides with $z$: 
$\tau_0=z$ (here restriction  ${\rm Re}(\beta)>0$ 
should be annulated). 
It provides convergence at the point $\tau_0$ for $0<\tau_0\leq z$. 

Using notations \eqref{chi}, \eqref{psi}, \eqref{min} and \eqref{tau} for $\chi$, $\psi_p$, ${\rm M}$ and $\tau_0$ we  can write the following conditions of convergence for integral \eqref{GINT}:
\begin{equation}\label{cond07}
	\left\{
	\begin{array}{ll}
		{\rm Re}(\beta )>0,\,\, {\rm M}>0,  & p<q; \\
		{\rm M}>0,\,\,	\tau_0\in \mathbb{R} \,\,\text{and}\,\,\\
		\left\{ \begin{array}{ll}
			0<\tau_0<z,\,\,{\rm Re}\left(\psi_p\right)>-1 & \mbox{if ${\rm Re}(\beta )>0$};\\
			\tau _0=z & \mbox{if ${\rm Re}\left(\beta +\psi_p\right)>0$},\end{array} \right.& p=q;\\
		{\rm Re}(\beta )>0, \,\,{\rm M}>0\,\,\text{and}\,\,\\
		\left\{ \begin{array}{ll}
			{\rm Re}\left(\alpha +\frac{g  }{r}\chi\right)>0 & \mbox{if $\frac{g}{r}<0$};\\
			{\rm Re}\left((-1)^{q-m-n}w^{-\frac{1}{r}} \right)\leq 0 & \mbox{if $\frac{g}{r}>0$},\end{array} \right.& p=q+1; \\
		{\rm Re}(\beta )>0, \,\,{\rm M}>0,\,\, \sqrt{(-1)^{q-m-n-1} w^{-\frac{1}{r}} }\in \mathbb{R},\,\\ \frac{g}{r}\in \mathbb{R},\,{\rm Re}\left(\frac{\alpha 
			r}{g}+\chi \right)<\frac{1}{2},& p=q+2; \\
		{\rm Re}(\beta )>0,\,\,{\rm M}>0,\,\,\text{and}\\
		\left\{ \begin{array}{ll}
			{\rm Re}\left(\alpha +\frac{g}{r} \chi \right){>}0 & \mbox{if $\frac{g}{r (q-p)}{>}0$};\\
			{\rm Re}\left((p-q) \left(-w^{\frac{1}{r}}\right)^{\frac{1}{q-p}}\right){\leq}0 & \mbox{if $\frac{g}{r (q-p)}{<}0$},\end{array} \right.& p>q+2, \\
	\end{array}
	\right.
\end{equation}

\subsection{ Support of differential constants. Generic formulas for fractional differentiation}\label{BF}

With the advent of computer algebra systems, such as Wolfram Mathematics, Maple, Matlab, etc., it became necessary to revise approaches to well-known functions. Computer systems demanded to provide correct numerical evaluations of analytical functions everywhere in complex plane, including on branch lines like $(-\infty,0)$ for $\sqrt{z}$. It stimulated developers of the Mathematica system to revise mathematical formulas, where the behavior of functions on branch cuts was accurately described not only theoretically, but supported by numerical evaluations in Mathematica. As a result, Wolfram Mathematics uses a simple axiom: "The argument of all complex numbers $z$ satisfies the inequality $-\pi<{\rm arg}(z)\leq\pi$."\, As a result functions are internally consistent with each to other and can be described  even on branch cut lines.
As a result \texttt{Mathematica} makes operations in the full complex plane with all functions involving so-called differential constants (like $\frac{\sqrt{z^2}}{z}$, $\log \left(z^2\right)-2 \log (z)$) and    piecewise construction including logarithmic situations. 
The question of what happens to the fractional integro-differentiation in such cases previously ignored in the literature. Operation \texttt{FractionalOrderD} fills this gap.

By definition in Mathematica we have
\begin{equation}\label{Const01}
	\frac{\sqrt{z^2}}{z}=\left\{
	\begin{array}{cc}
		1, & -\frac{\pi }{2}<\arg (z)\leq \frac{\pi }{2}; \\
		-1, & \text{in other cases} \\
	\end{array}
	\right.
\end{equation}
and
\begin{equation}\label{Const02}
	\log \left(z^2\right)-2 \log (z)=\left\{
	\begin{array}{cc}
		0, & -\frac{\pi }{2}<\arg (z)\leq \frac{\pi }{2}; \\
		2 i \pi,  & -\pi <\arg (z)\leq -\frac{\pi }{2}; \\
		-2 i \pi,  & \text{in other cases}. \\
	\end{array}
	\right.
\end{equation}

The built-in  Mathematica derivative operator \texttt{D} misses discontinuities at branch cuts:
\begin{equation}\label{Const03}
	\text{{\bf D}}\left[\frac{\sqrt{z^2}}{z},z\right]=0
	\qquad \text{and} \qquad
	\text{{\bf D}}\left[\log \left(z^2\right)-2 \log (z),z\right]=0.
\end{equation}

\texttt{FractionalOrderD} provides the following results, which saves components $\sqrt{z^2}$ and $\log(z^2){=}\log(-iz){+}\log(iz)$:
\begin{equation}\label{Const04}
	\frac{d^\alpha}{dz^\alpha} \frac{\sqrt{z^2}}{z} =\frac{\sqrt{z^2} z^{-\alpha -1}}{\Gamma (1-\alpha )}
\end{equation}
and 
\begin{equation}\label{Const05}
	\frac{d^\alpha}{dz^\alpha}[\log \left(z^2\right)-2 \log (z)]=\frac{z^{-\alpha } (\log (-i z)+\log (i z)-2 \log (z))}{\Gamma (1-\alpha )}.
\end{equation}

Classical (order $1$ or $2$, etc) differentiation and integration can operate with abstract functions $f(z),g(z),h(z)$ and some their constructions (product, composition, ratio, inverse function, series, etc.). Similar formulas can be derived for fractional integro-differentiation. Many of such formulas we already have in \texttt{ResourceFunction["FractionalOrderD"]}. For example, we have generic formulas for 
$$
f^{(-1)}(z),\qquad f(z)\cdot g(z),\qquad f(z)\cdot g(z)\cdot h(z),\qquad 
\frac{1}{f(z)},\qquad
\frac{f(z)}{g(z)},\
$$
$$
(f(z))^a,\qquad (f(z))^{g(z)},\qquad f(z)^{g(z)^{h(z)}},\qquad (a z)^b \log ^c(d\cdot f(z)),
$$
$$
c^{b f(z)^a+d},\qquad c^{b f\left(z^a\right)+d},\qquad c^{b f\left(a^z\right)+d},\qquad f(g(z)),\qquad f(g(h(z))),\qquad\text{etc.}
$$
Here $f^{(-1)}(z)$ is an inverse for $f(z)$ function. 
Below we give one rather simple example for product $f(z)\cdot g(z)$: 
\begin{equation}\label{Mult01}
	\frac{d^\alpha }{dz^\alpha}(f(z)g(z))= \left\{ \begin{array}{ll}
		\sum\limits_{k=0}^{\alpha } \binom{\alpha }{k} f^{(\alpha-k)}(z) g^{(k)}(z)  & \mbox{if $z\in\mathbb{Z},\alpha\geq 0$};\\
		\sum\limits_{k=0}^{\alpha } \binom{\alpha }{k}  \frac{d^{\alpha
				-k}f(z)}{d z^{\alpha -k}}\frac{d^{k}g(z)}{d z^{k}} & \mbox{in other cases}.\end{array} \right.
\end{equation} 

Of course, most formulas of fractional integro-differentiation of abstract functions are very complicated and large. 
But in the case of specific formulas, we get simpler expressions. For example, 
\begin{equation}\label{Exp01}
	\frac{d^\alpha }{dz^\alpha}e^{z^2}=\frac{\sqrt{\pi } 2^{\alpha } z^{-\alpha }}{\Gamma\left(1-\frac{\alpha+1}{2}\right) \Gamma\left(1-\frac{\alpha
		}{2} \right) } \, _2{F}_2\left(\frac{1}{2},1; 1-\frac{\alpha+1}{2},1-\frac{\alpha
	}{2};z^2\right).
\end{equation}

In recent years Wolfram Language has added some new generic functions 
\begin{itemize}
	\item 10 Heun functions,
	\item 4 Lam\'e functions,
	\item 8 Carlson elliptic integrals,
	\item Fox $H$--function,
	\item  4 Coulomb functions.
\end{itemize}
For these functions   corresponding formulas of fractional order integro-differentiation can be obtained. 
For example,
we consider Carlson's elliptic integral (see \cite{Carlson})
\begin{equation}\label{Carlson01}
	R_D(x,y,z)=\frac{3}{2}\int\limits_0^\infty (t+x)^{-1/2}(t+y)^{-1/2}(t+z)^{-3/2}dt,\quad x>0,\quad y>0, \quad z>0 
\end{equation} 
and  Heun $G$-function $H(a,b,c,d,p,q,z)$ (see \cite{Heun}). 

Heun $G$-function $H(a,b,c,d,p,q,z)$ specializes to $\,_2F_1(c,d,p;z)$ if $b=a\cdot c\cdot d$ and $q=c+d-p+1$ or $a=1$ and $b=c\cdot d$.

We obtain
$$
\frac{d^\alpha }{dz^\alpha}R_D\left(\frac{1}{2},z,1\right)=
$$
\begin{equation}\label{Carlson02}
	\begin{cases}
		\frac{3 (-1)^{\alpha }}{2 \alpha +3} \left(\frac{1}{2}\right)_{\alpha } F_1\left(\alpha +\frac{3}{2};\frac{1}{2},1;\alpha
		+\frac{5}{2};\frac{1}{2},1-z\right), & \alpha \in \mathbb{Z}, \alpha \geq 0; \\
		\frac{3 \pi  z^{-\alpha }}{4 \Gamma (1-\alpha )} \sum\limits_{k=0}^{\infty } \frac{\,
			_2F_1\left(\frac{1}{2},k+\frac{3}{2};2;\frac{1}{2}\right) \left(\frac{3}{2}\right)_{k}
			z^{k}}{(1-\alpha )_{k}}\\
		-\frac{3 \sqrt{\pi } z^{\frac{1}{2}-\alpha }}{2 \Gamma \left(\frac{3}{2}-\alpha \right)} \sum\limits_{k=0}^{\infty } \frac{(k+1)! \, _2F_1\left(\frac{1}{2},k+2;2;\frac{1}{2}\right)
			z^{k}}{\left(\frac{3}{2}-\alpha \right)_{k}}, &
		\text{in other cases}
	\end{cases}
\end{equation} 
and
\begin{equation}\label{Heun01}
	\frac{d^\alpha }{dz^\alpha}H(a,b,c,d,p,q,z)=\sum _{k=0}^{\infty } \frac{k! c_k z^{k-\alpha }}{\Gamma (k-\alpha +1)},
\end{equation} 
where $| z| <\min (1,| a| )$, $ c_0=1$, $c_1=\frac{b}{a p}$, $c_j=-\frac{c_{j-2} P_{j-2}+c_{j-1} Q_{j-1}}{R_j}$, $P_j=(c+j) (d+j)$, $Q_j=-j (a (j+p+q-1)+c+d+j-q)-b$,
$R_j=a j (j+p-1)$.

Heun $G$-function, in general, is not hypergeometric type function and its series representation has coefficients $c_k=\frac{k! c_k}{\Gamma (k-\alpha +1)}$, $\alpha=0$, which satisfies three-terms recurrent relations, which were described in \texttt{ConditionalExpression}.

Currently we  have   general rules for:
$a^{g(z)}$, $\left(a^{b f(z)}\right)^c$, $c^{b a^{f(z)}+d}$, $c^{b f\left(a^z\right)+d}$, $c^{b f(z)^a+d}$, $c^{b f\left(z^a\right)}$, $e^{g(z)}$, $e^{g(h(z))}$, $f\left(a^z\right)$, $(f\left(a^{b z}\right))^c$, $f\left(a^{b z^c}\right)$, $\frac{1}{f(z)}$, $z f(z)$, $z^2 f(z)$, $(f(z))^a$, $(f(z))^{g(z)}$, $(f(z))^{(g(z))^{h(z)}}$, $f\left(z^a\right)$, $f(g(z))$, $f(g(h(z)))$, $\frac{f(z)}{g(z)}$, $f(z)\cdot g(z)$, $f(z)\cdot g(z)\cdot h(z)$, $f^{(-1)}(z)$. 

\section{Conclusion}

Having considered various approaches to introducing an arbitrary power of a differential operator $\frac{d}{dx}$, we come to the conclusion that these approaches are not so different. 
So, when calculating various fractional derivatives of a power function $x^p$, we almost always get the same result. 
This makes the considered derivatives coincide on the class of analytic functions. On the other hand, if there are additional parameters in the fractional derivative, then the result of applying such an operator to a power function will depend on these parameters. In some applied problems, this makes sense. Another reason for introducing different definitions of fractional derivatives is to apply regularizations of divergent integrals. Wolfram Mathematica uses the Hadamard regularization of divergent integrals. In addition, Mathematics systematizes formulas for fractional derivatives of general form and for general functions.

Finally, returning to main question: "How maximally naturally make extensions of differentiation $\frac{d^n}{dx^n}$ from natural order $n=1,2,3,...$ to arbitrary symbolical (or complex) order $\alpha$"\, we can say the following: "Like in the case of extension factorial $n!$ from natural $n$ to Gamma function $\Gamma(\alpha+1)$  ($\alpha=n$) we do not have unique solution in such extension (corresponding details for $n!$ see in  \cite{Marichev}, \cite{Marichev01} page 43.). As result we see in literature the numerous approaches to definitions of fractional order differentiation by variable $x$. Corresponding inversions define fractional order integrations, which in majority cases include additional point a at the end of integral from a to x, which demands additional restrictions for convergence. Existence of this point distracts attention and it is naturally to use $a=0$ as the end of integral from $0$ to $x$ and allow $x$ to be negative or complex. Such construction leads to Riemann-Liouville left-sided fractional integral with beginning at $0$ and Riemann-Liouville fractional derivative described at formula \eqref{FrW}. After involving Hadamard's concept of the "finite part"\, we avoid influence of other "aside" point ($x=0$) on procedure of fractional integro-differentiation in our main point $x$, which effectively works in all complex $x$-plane for any analytical functions in their regular or branch points, including logarithmic cases ($\log(x)$, $\frac{1}{x}$,  $x^c$, $(x^n)^{1/m}$, etc.). It allowed to form described above Riemann-Liouville-Hadamard fractional order integro-differentiation, which looks as maximally natural extension of  integer order differentiation and its inversion (integer order repeatable indefinite integration without arbitrary polynomials, that arises). This statement follows from comparison of different approaches to fractional order integro-differentiation between themself and how they realized on the basic simplest functions $z^\lambda$, $e^z$, etc. which is widely used in Taylor and Fourier series."      

Additional information about Fractional Differentiation one can find in the talk of Oleg Marichev, Paco Jain "New in Fractional Differentiation"\, at Wolfram virtual Technology Conference (October 12-15, 2021) \cite{Ev01}, which was extended to the next talk of Oleg Marichev "Compositional Structure of Classical Integral Transforms" at the Wolfram Technology Conference (October 18-21, 2022) \cite{Ev02}. Later these results got development at the talk of Marichev O., Shishkina E. \cite{MarSh} and set of Internet presentations of Oleg Marichev \cite{Com01,Com02,Com03}.


\end{document}